\documentclass[12pt]{article}
\usepackage{amssymb,latexsym, amsmath, enumerate}
\usepackage{color}

\def\U{{\mathrm{U}}}
\def\X{{\mathrm{X}}}
\def\Y{{\mathrm{Y}}}
\def\R{\mathbb R}
\def\N{\mathbb N}
\def\C{\mathbb C}
\def\H{\mathbb H}

\newtheorem{theorem}{Theorem}
\newtheorem{lemma}{Lemma}
\newtheorem{remark}{Remark}

\newtheorem{corollary}{Corollary}
\newtheorem{proposition}{Proposition}

\title{Revisiting the heat kernel on isotropic and nonisotropic Heisenberg groups \\Dedicated to Professor Jia-Rong Yu on his 100th birthday}
\author{Hong-Quan Li, Ye Zhang}
\date{}
\begin{document}

\renewcommand{\theequation}{\thesection.\arabic{equation}}
\setcounter{equation}{0} \maketitle

\vspace{-1.0cm}

\bigskip

{\bf Abstract.} The aim of this paper is threefold. First, we obtain the precise bounds for the heat kernel on isotropic Heisenberg groups by using well-known results in the three dimensional case. Second, we study the asymptotic estimates at infinity for the heat kernel on nonisotropic Heisenberg groups. As a consequence, we give uniform upper and lower estimates of the heat kernel, and complete its short-time behavior obtained by Beals-Gaveau-Greiner. Third, we complete the results obtained in \cite{Li12} about the heat kernel of Grushin operators.

\medskip

{\bf Mathematics Subject Classification (2000):} {\bf  58J37, 35B40, 35H10, 35B45, 35K05, 43A80, 58J35, 43A85}

\medskip

{\bf Key words and phrases:}  Heat kernel;
Asymptotic behavior; Heisenberg groups; Precise estimates; Grushin operators

\medskip

\renewcommand{\theequation}{\thesection.\arabic{equation}}
\section{Introduction}
\setcounter{equation}{0}

The small-time asymptotics of the heat kernel in the setting of Lie groups as well as hypoelliptic operators have attracted considerable attention during the past three decades, see for example \cite{L87,BA88,BGG00,BB09,Li10,Li12,B12,BBN12,BBN16}
and references therein. They involve some geometric properties of the underlying space, such as the Carnot-Carath\'eodory distance and cut locus, etc.

In this paper, we consider related problems on a stratified Lie group $G$ endowed with the canonical sub-Laplacian $\Delta$ and the group of dilations $\delta_r$ ($r > 0$). Let $Q$ denote its homogeneous dimension. It is well-known that the heat kernel $p_h$ ($h > 0$), that is the convolution kernel of
$e^{h \Delta}$, satisfies the following scaling property (cf. for example \cite{FS82})
\begin{align} \label{hqp}
0 <  p_h(g) = h^{-\frac{Q}{2}} p_1(\delta_{\frac{1}{\sqrt{h}}}(g)), \qquad  \forall h > 0, \  g \in G.
\end{align}
Let $o$ denote the neutral element of $G$. The semigroup property implies that
\begin{align} \label{hup}
p_h(g) \leq p_h(o) = C(G, \Delta) h^{-\frac{Q}{2}}, \quad \forall g \in G, \ h > 0,
\end{align}
for a constant $C(G, \Delta)$ depending on $G$ and $\Delta$.
Now the problem of the short-time behaviour of the heat kernel can be reformulated by, for a given point $g \neq o$ (by \eqref{hup}, the case $o$ is trivial but abnormal), to find $\Theta_h(g) > 0$ (usually of type $C(g) h^{-\sigma(g)}$) such that
\begin{align} \label{STA}
p_h(g) = \Theta_h(g)  e^{-\frac{d^2(g)}{4 h}} \left( 1 + o(1) \right), \qquad \mbox{as $h \longrightarrow 0^+$,}
\end{align}
where $d(g) = d(g,o)$ denotes the Carnot-Carath\'eodory distance (associated to $\Delta$) between $g$ and $o$. Here and in the sequel we use the notation $f = o(w)$ if $\lim \frac{f}{w} = 0$. Recall that we have
\begin{align} \label{DS}
d(\delta_r(g)) = r d(g), \quad \forall r > 0, \  g \in G.
\end{align}

In fact, \eqref{STA} is a refinement of Varadhan's formula, which is valid in a general frame (cf. \cite{L87}, \cite{AH05} and references therein), as follows
\begin{align}\label{Varadhan}
\lim_{h \longrightarrow 0^+} 4 h \, \ln{p_h(g)} = - d^2(g), \quad \forall g \in G.
\end{align}

To our knowledge, many works on \eqref{STA} can be found in the setting of stratified Lie groups, but a complete and precise calculation has been obtained only for H-type groups, see \cite{G77} for the three dimensional Heisenberg group, \cite{BGG00} or \cite{Li07} for isotropic Heisenberg groups, and \cite{Li10} for other H-type groups. We also recall that some partial results have been obtained for nonisotropic Heisenberg groups by Beals, Gaveau and Greiner in \cite{BGG00} (see also \cite{BGG99}). One goal of this article is to complete their results.

We point out that in contrast to the elliptic case, \eqref{STA} cannot be uniform in a neighbourhood of the neutral element, see for example \cite{G77} or the results below.
Although there exists scaling invariance property, cf. \eqref{hqp}, the small-time asymptotic behavior of the heat kernel does not directly imply Gaussian estimates of the heat kernel which play an important role in analysis.

Certainly, there exist some effective methods to obtain off-diagonal Gaussian upper bounds (that is the Gaussian upper bounds of $p_h(g)$ here) from on-diagonal upper bounds  (namely, the upper bounds of $p_h(o)$ in our situation). These methods are valid in some much more general situations than stratified groups, including but not limited to the following: (1) logarithmic Sobolev inequalities; (2) Nash type inequalities; (3) Sobolev inequalities; (4) the finite propagation speed property for the associated wave equation; (5) the Phragm\'en-Lindel\"of theorem. See for example \cite{DP89, C93, S96, S04, CS08, BCS15} and the references therein. In particular, the following upper estimates  %which are sharp due to Molchanov's result on $Q$-dimensional sphere
are well-known:
\begin{align} \label{HU1}
p_h(g) \lesssim h^{-\frac{Q}{2}} \left( 1 + \frac{d(g)}{\sqrt{h}} \right)^{Q - 1} e^{-\frac{d^2(g)}{4 h}}, \quad \forall h > 0, \ g \in G,
\end{align}
where the power $Q - 1$ in the second term is optimal due to Molchanov's result on $Q$-dimensional sphere (see \cite[Example 3.1, p. 23]{M75}).
Here and in the sequel we use the notation $f \lesssim w$ (resp. $f = O(w)$) if there exists a universal constant $C > 0$ (maybe depending on $(G, \Delta)$) such that $f \leq C w$ (resp. $|f| \leq C \, w$).

Classical lower estimates for the heat kernel, that is, for any $0 < \varpi < 4$, there exists a constant $C(\varpi) > 0$ such that
\begin{align} \label{HL1}
p_h(g) \geq C(\varpi) h^{-\frac{Q}{2}} e^{-\frac{d^2(g)}{(4 - \varpi) h}}, \quad \forall h > 0, \ g \in G,
\end{align}
as well as upper bounds for its derivatives can be found in \cite{V90, VSC92} and \cite{CS08}. Those classical estimates and their improvements \eqref{HU1} are far enough for studying problems such as Hardy-Littlewood-Sobolev inequalities, Riesz transforms, etc.

Recall that we have the inequality \eqref{hup}, and also notice that \eqref{HU1} follows from some general approach
and reveals little
geometric information  of the underlying groups. Hence it is not hard to understand why some results, related to dispersive and Strichartz estimates for wave or Schr\"odinger operators as well as spectral multipliers obtained on step-2 stratified groups (see \cite{MS94, DH05, MS15, MM16, BFKG16} and the references therein), could not be explained by \eqref{HU1}. However, asymptotic estimates of the heat kernel obtained on isotropic Heisenberg groups (cf. \cite{Li07}) and on H-type groups (cf. \cite{Li10}) shed light on those subjects. Of course, further researches are necessary to provide new proofs for results cited above via heat kernel.

We recall the gradient bound of the heat semigroup on the simplest non-trivial stratified Lie group, the three dimensional Heisenberg group $\H(1, 1)$: (see for example \cite{B94, B97, Le00, ABCGMRS00, Li06, BBBC08, Q10, LP10, Le11, Q11, Q13, BGL14, EGS18} and the references therein for the background and its extensive applications)
\begin{align} \label{HGE}
|\nabla e^{h \Delta} f|(g) \lesssim e^{h \Delta} \left( |\nabla f| \right)(g), \quad \forall h > 0, g, \ f \in C_0^{\infty},
\end{align}
where $\nabla$ denotes the horizontal gradient. To our best knowledge, there exist three different approaches to prove \eqref{HGE}. One is provided by the first author in \cite{Li06} and two others are given by Bakry-Baudoin-Bonnefont-Chafa\"{\i} in \cite{BBBC08}. All these methods are based on precise bounds of the heat kernel
and its gradient on $\H(1, 1)$ (see \cite{HM89}, \cite{BGG00}, \cite{Li06} or \cite{Li07}), that is
\begin{align}
p_h^{\H(1, 1)}(z, t) \sim  \frac{h^{-2} e^{-\frac{d^2(g)}{4 h}}}{\sqrt{1 + \frac{|z| d(g)}{h}}}, \quad \forall h > 0, & \ g = (z, t) \in \H(1, 1) \cong \C \times \R, \label{HEB1} \\
|\nabla \ln{p_h^{\H(1, 1)}(g)}| \lesssim & \frac{d(g)}{h}, \quad \forall h > 0, g \ \in \H(1,1),  \label{HEB2}
\end{align}
where we adopt the notation $f \sim w$ if $f \lesssim w$ and $w \lesssim f$.

Notice that \eqref{HGE} has been generalized
to H-type groups in \cite{HL10} (resp. \cite{E10}) via the same method of the first author in \cite{Li06} (resp. of Bakry-Baudoin-Bonnefont-Chafa\"{\i} in \cite{BBBC08}). Recall the key ingredients are also precise bounds for the heat kernel and its gradient obtained in \cite{Li10} and  \cite{E09}, respectively. For other applications of precise bounds of the heat kernel, see for example \cite{HZ10} and \cite{BCH18}.

By the scaling invariance property, cf. \eqref{hqp}, estimates \eqref{HU1} and \eqref{HL1}, we can get small-time asymptotic and precise estimates immediately if the asymptotic behaviour at infinity for the heat kernel at time $1$, that is, $p_1(g)$ is known. However, the subject becomes much
more difficult and the only known case is H-type groups, see \cite{Li10}, also \cite{HM89} for the three dimensional Heisenberg group and \cite{Li07} for general isotropic Heisenberg groups.

The main purpose in this article is to obtain uniform asymptotic bounds at infinity for $p_1(g)$ on more complicated groups, the nonisotropic Heisenberg groups, and to shed additional light on the heat kernel on step-$2$ stratified groups.

\medskip

To begin with, let us recall some definitions and basic facts.

\subsection{Preliminaries on isotropic and non-isotropic Heisenberg groups}

Suppose that
\begin{align*}
l \in \N^+, \  \mathcal{K} = (k_1, \cdots, k_l) \in \left( \N^+ \right)^l,  \  n = \sum_{i = 1}^l k_i,
\end{align*}
and
\begin{align*}
\mathcal{A} = (a_1, \cdots, a_l)  \mbox{ with } 0 <  a_1 < \cdots < a_l = 1.
\end{align*}
Recall that (see \cite{BGG00} or \cite{LP03}) the general Heisenberg group with parameter $(\mathcal{K}, \mathcal{A})$,  $$\H(\mathcal{K}, \mathcal{A}) = \left( \prod_{i = 1}^l \C^{k_i} \right) \times \R,$$ is a stratified group
with the group multiplication
\begin{eqnarray*}
(z, t) \cdot (z', t') = (z + z', t + t' + 2 \sum_{i = 1}^l a_i \Im
\langle z_i, z_i' \rangle),
\end{eqnarray*}
where
$$z = (z_1, \cdots, z_l),  \quad z' =
(z_1', \cdots, z_l') \in \prod_{i = 1}^l \C^{k_i},$$
$\Im \omega$ denotes the imaginary part of a complex number $\omega$ and $\langle \cdot, \cdot \rangle$ denotes the usual complex inner product. When $l = 1$, $\H(n, 1)$ is often called the Heisenberg group or the isotropic Heisenberg group of real dimension $2 n + 1$. When $l \geq 2$, $\H(\mathcal{K}, \mathcal{A})$ is often called the non-isotropic Heisenberg group of real dimension $2 n + 1$. Recall that the dilation $\delta_r$ ($r > 0$) is given by
\begin{align*}
\delta_{r}(z, t) = (r z, r^2 t), \quad \forall (z, t) \in \H(\mathcal{K}, \mathcal{A}) = \C^n \times \R.
\end{align*}

For $1 \leq i \leq l$, set
$$z_i = (z_{i, 1}, \cdots, z_{i, k_i}) = (x_{i, 1} + \imath y_{i, 1}, \cdots, x_{i, k_i} + \imath y_{i, k_i}) \in \C^{k_i} $$
with
$$x_{i, j}, y_{i, j} \in \R, \quad  1 \leq j \leq k_i, \quad |z_i| = \sqrt{\langle z_i, z_i \rangle}, \quad |z| = \sqrt{\sum_{i = 1}^l |z_i|^2}.$$
Let
\begin{eqnarray*}
\X_{i, j} = \frac{\partial}{\partial x_{i, j}} + 2 a_i y_{i, j} \frac{\partial}{\partial t},  \quad \Y_{i, j} = \frac{\partial}{\partial
y_{i, j}} - 2 a_i x_{i, j} \frac{\partial}{\partial t}, \qquad 1 \leq i \leq l, 1 \leq j \leq k_i,
\end{eqnarray*}
be the left invariant vector fields on $\H(\mathcal{K}, \mathcal{A})$. The associated canonical sub-Laplacian and
horizonal gradient are given respectively by
\begin{eqnarray*}
\Delta_{\H(\mathcal{K},\mathcal{A})} = \sum_{i = 1}^l \sum_{j = 1}^{k_i} (\X_{i, j}^2 + \Y_{i, j}^2), \qquad \nabla_{\H(\mathcal{K},\mathcal{A})} = (\X_{1, 1}, \Y_{1, 1}, \cdots, \X_{l, k_l}, \Y_{l, k_l}).
\end{eqnarray*}

Let $d(g) = d_{\H(\mathcal{K},\mathcal{A})}(g)$ denote the Carnot-Carath\'eodory distance between  $g$ and the neutral element $o$,
$p_h^{\H(\mathcal{K},\mathcal{A})}$ ($h
> 0$) the heat kernel, namely the convolution kernel of $e^{h \Delta_{\H(\mathcal{K},\mathcal{A})}}$, and we set $p^{\H(\mathcal{K},\mathcal{A})} = p^{\H(\mathcal{K},\mathcal{A})}_1$.
Recall that for all $h > 0$ and $(z, t) \in \H(\mathcal{K}, \mathcal{A})$, we have \begin{eqnarray} \label{sp}
d\left(\frac{z}{\sqrt{h}},\frac{t}{h}\right)=\frac{1}{\sqrt{h}}d(z,t),  \quad p^{\H(\mathcal{K},\mathcal{A})}_h(z, t) = \frac{1}{h^{n+1}} p^{\H(\mathcal{K},\mathcal{A})}\left(\frac{z}{\sqrt{h}}, \frac{t}{h}\right),
\end{eqnarray}
and
\begin{align}
p^{\H(\mathcal{K},\mathcal{A})}(z, t) &= \frac{1}{2 (4\pi)^{n+1}} \nonumber \\
& \times \int_{\R}\prod_{j=1}^l \left( \frac{a_j \lambda}{\sinh{(a_j \lambda)}} \right)^{k_j} \exp{\left( \frac{1}{4} \left(\imath \lambda t - \sum_{j=1}^l |z_j|^2 a_j \lambda \coth{(a_j \lambda)}\right) \right)} \, d\lambda \label{1}
\end{align}
by \cite{BGG00} or \cite{LP03}.

\medskip

Then we provide a new approach, which is of independent interest, to yield the following.\\

\subsection{Precise estimates for the heat kernel on isotropic Heisenberg groups $\H(n, 1)$}

\begin{theorem} \label{TC}
For $n \ge 2$, we have
{\em \begin{align} \label{IHE}
p_h^{\H(n, 1)}(z, t) \sim h^{- n - 1} e^{-\frac{d^2(z, t)}{4 h}} \left( 1 + \frac{d(z, t)}{\sqrt{h}} \right)^{2 (n - 1)} \left( 1 + \frac{|z| d(z, t)}{h} \right)^{\frac{1}{2} - n},
\end{align}}
for all $h > 0$ and $(z, t) \in \H(n, 1) \cong \C^n \times \R$.
\end{theorem}

We emphasize that the implicit constants in \eqref{IHE} depend on $n$. Remark that the estimate is valid for all $n \ge 1$:
the upper bound was obtained under an equivalent form by Beals-Gaveau-Greiner (see \cite[Theorem 4.6]{BGG00}), and \eqref{IHE} was obtained in \cite{Li07}. Moreover, the uniform asymptotic behaviour at infinity of $p^{\H(n, 1)}$ for any $n \ge 1$ has been acquired by the first author in \cite{Li07}, by improving the method of Beals-Gaveau-Greiner in \cite{BGG00}. It is
also obtained in \cite{Li12} by another approach which is more powerful since $p^{\H(n, 1)}$ can be considered as a special case of the heat kernel of the  Grushin operator $\Delta_G$ on $\R^{2 n} \times \R$, where
\begin{align*}
\Delta_G = \sum_{i = 1}^{2 n} \left( U_i^2 + V_i^2 \right), \quad \mbox{with } U_i = \frac{\partial}{\partial x_i}, \  V_i = x_i \frac{\partial}{\partial t}.
\end{align*}
Strictly speaking, the asymptotic behaviour of $p^{\H(1, 1)}$ cannot be deduced directly from the results in \cite{Li12}. However, we will see in the following that an improvement of the method in \cite{Li12} is sufficient for our desired result. Recall that the three dimensional Heisenberg group case has been proved by Hueber-M\"uller in \cite{HM89}.

Now we consider our third proof of the above result. First we observe from the expression of the heat kernel  given above (cf. \eqref{1}) and an elementary property of Fourier transform (see \eqref{Fourier} below for more details) that for any $n \ge 2$ and $(z,t) \in \H(n,1)$, we have
\begin{align}
p^{\H(n, 1)}(z, t) &=  \left( p^{\H(1, 1)}(|z|, \cdot) \ast p^{\H(n - 1, 1)}(o, \cdot) \right)(t) \label{convolution} \\
&=  \left( p^{\H(1, 1)}(|z|, \cdot) \ast \underbrace{p^{\H(1, 1)}(0, \cdot) \ast \cdots \ast p^{\H(1, 1)}(0, \cdot)}_{n - 1 \, \mbox{times}} \right)(t) \nonumber
\end{align}
where the convolution $\ast$ is taken in $\R$.

The last formula can be regarded as a reinforced version of the following identity:
\begin{align}\label{convolution2}
p^{\H(n, 1)}(z, t) = \left( p^{\H(1, 1)}(z_1, \cdot) \ast p^{\H(1, 1)}(z_2, \cdot) \ast \cdots \ast p^{\H(1, 1)}(z_n, \cdot) \right)(t),
\end{align}
which has been used in \cite[p. 371]{C96}.
Then our proof is based on the above observation, some results of Gaveau in \cite{G77} and the estimate \eqref{HEB1}.

\medskip

In order to state our main results,
we fix some notations.

\subsection{Notations}\label{reference}

Let $w$ be a non-negative real-valued function. From now on, the notation $f = O(w)$ (resp. $f \lesssim w$ if $f$ is also real-valued function) means that there exists a universal constant $C > 0$, maybe depending on $(\mathcal{K}, \mathcal{A})$, such that $|f| \leq C w$ (resp. $f \leq C w$). Recall that we use the notation $f \sim w$ if $f \lesssim w$ and $w \lesssim f$. Similarly, $f = O_{\vartheta}(w)$ (resp. $f \lesssim_{\vartheta} w$ for real-valued function $f$), with a parameter $\vartheta$, means that there exists a constant $C(\vartheta) > 0$, depending on $\vartheta$ and maybe depending on $(\mathcal{K}, \mathcal{A})$, such that $|f| \leq C(\vartheta) w$ (resp. $f \leq C(\vartheta) w$). Furthermore, we use the notation $f \sim_{\vartheta} w$ if $f \lesssim_{\vartheta} w$ and $w \lesssim_{\vartheta} f$. We remark that starting from Section \ref{st1n}, $|f|$ often denotes the magnitude of the complex-valued function $f$.

\medskip

Then we recall the Carnot-Carath\'{e}odory distance on $\H(\mathcal{K}, \mathcal{A})$ obtained in \cite[\S 3]{BGG00}.
In fact, this part is not necessary for our main goal, namely Theorems \ref{t1}-\ref{t3}
in Subsection \ref{mainre} below and we include it for the sake of clarity
(see the discussion in the end of Subsection \ref{mainre}).

Set
\begin{align} \label{mu}
\mu(\omega) = \frac{\omega}{{\sin^2{\omega}}}-\cot{\omega} = \frac{2 \omega - \sin{2 \omega}}{2 \sin^2{\omega}}: \quad (-\pi, \pi) \longrightarrow \R,
\end{align}
which is an odd function, and a monotonely increasing diffeomorphism between $(-\pi,\pi)$ and $\R$. Given $(z,t)$ such that $z_l \ne 0$, there exists exactly one $-\pi < \theta = \theta(z, t) < \pi$ such that
\begin{eqnarray} \label{dd1}
t = \sum_{j=1}^l a_j \mu(a_j \theta) |z_j|^2,
\end{eqnarray}
and we have
\begin{eqnarray} \label{dd2}
d^2(z,t) = \sum_{j=1}^l \left( \frac{a_j \theta}{\sin(a_j \theta)} \right)^2 |z_j|^2
= \theta \, ( t + \sum_{j=1}^l a_j \cot{(a_j\theta)} \, |z_j|^2 ).
\end{eqnarray}

We have the following alternatives:
1. If $|t| < \sum\limits_{j=1}^{l - 1}
a_j \mu(a_j \pi) |z_j|^2$,  then the equation \eqref{dd1} has a unique solution $-\pi < \theta < \pi$ and  \eqref{dd2} remains valid.
2. If $|t| \geq \sum\limits_{j=1}^{l - 1}
a_j \mu(a_j \pi) |z_j|^2$,  then we have
\begin{eqnarray*}
d^2(z, t) =
\pi \left( |t| + \sum_{j=1}^{l-1} a_j \cot{(a_j \pi)} |z_j|^2 \right).
\end{eqnarray*}

Using the equivalence between the Carnot-Carath\'eodory distance and a homogeneous norm on stratified groups (see for example \cite{VSC92}), or by a
direct calculation, we have
\begin{align} \label{l1}
d^2(z, t) \sim |z|^2 + |t|, \qquad \forall (z, t) \in \H(\mathcal{K}, \mathcal{A}).
\end{align}

\medskip

Recall that $p^{\H(\mathcal{K},\mathcal{A})}  \in C^{\infty}(\H(\mathcal{K} , \mathcal{A}))$ and note that
\begin{align*}
p^{\H(\mathcal{K},\mathcal{A})}(z, t) = p^{\H(\mathcal{K},\mathcal{A})}(z, -t),
\quad d(z, t) = d(z, -t), \quad \forall (z, t) \in \H(\mathcal{K}, \mathcal{A}).
\end{align*}
Hence, without loss of generality, we may assume in the sequel that $t \geq 0$ and $z_l \ne 0$. Under our assumption, we get $0 \leq \theta < \pi$.

We define for $\lambda \in \C$ with $- \pi < \Im \lambda < \pi$,
\begin{eqnarray}\label{defvarphi}
\varphi(\lambda) = \varphi(\mathcal{K}, \mathcal{A}; z, t; \lambda)= \frac{1}{4} \left( \imath \lambda t -\sum_{j=1}^l |z_j|^2 a_j \lambda \coth{(a_j \lambda)} \right), \end{eqnarray}
and
\begin{align} \label{phi}
\Phi(s)=\Phi(\mathcal{K}, \mathcal{A}; z, t; s)=\varphi(s + \imath \theta) - \varphi(\imath \theta), \qquad s \in \R,
\end{align}
where $\theta$ is given by \eqref{dd1}.
Note that
\begin{align}
& \varphi(\imath \theta) = - \frac{d^2(z, t)}{4}, \quad  \varphi'(\imath \theta) = 0 = \Phi(0) = \Phi'(0) , \label{de1} \\
& \Phi^{\prime \prime}(0) =- \sum_{j=1}^l \frac{ a_j^2 |z_j|^2 }{ 2 \sin^2 ( a_j\theta ) }(1 - a_j \theta \cot( a_j \theta ) ) < 0. \label{dep}
\end{align}

In the case that $\theta$ is close to $\pi$, we set
\begin{align}\label{epsilon}
\epsilon = \pi - \theta.
\end{align}
We define holomorphic functions $G_i$ ($1 \leq i \leq 3$) and $G$ (with natural domain of definition) as follows:
\begin{eqnarray*}
G_1(\xi)=\frac{\xi}{\sin{\xi}} \cos{\xi}, \quad G_2(\xi)=\pi \left( \frac{\cos{\xi}}{\sin{\xi}} - \frac{1}{\xi} \right) = \pi \frac{\xi \cos{\xi} - \sin{\xi}}{\xi \sin{\xi}},
\end{eqnarray*}
and
\begin{align} \label{Ge}
G_3(\xi) &= G_3(\mathcal{K},\mathcal{A};z;\xi) = -\sum_{j=1}^{l-1} \frac{|z_j|^2}{4} G_1(a_j(\pi-\xi)), \\
\label{Ge2}
 G(\xi) &= G(\mathcal{K},\mathcal{A};z;\xi) =  \ G_3(\xi) + \frac{|z_l|^2}{4} (G_2(\xi) - G_1(\xi)).
\end{align}
It follows that
\begin{eqnarray} \label{v1p}
\varphi(\imath (\pi - \xi)) = -\frac{t}{4} (\pi - \xi) + G(\xi) + \frac{\pi}{4} \frac{1}{\xi} |z_l|^2,
\end{eqnarray}
with $\frac{\pi}{4} \frac{1}{\xi} |z_l|^2$ the bad part of the phase. We remark that
\begin{align}
G_2^{''}(0) = 0, \  G_1^{''}(0) = -\frac{2}{3}, \ G_1'(\xi) = -\mu(\xi), \  G_1^{''}(\xi) = -\mu'(\xi) < 0, \  \forall 0 \leq \xi < \pi.
\end{align}

We deduce from \eqref{v1p} and $\varphi'(\imath \theta) = 0$ that
\begin{align*}
\frac{t}{4}  =  - G'(\epsilon) + \frac{\pi}{4} \frac{|z_l|^2}{\epsilon^2}.
\end{align*}
Combining it with \eqref{v1p}, we get
\begin{align}
\Phi(\epsilon \tau) &= \varphi\left(\imath (\pi - \epsilon (1 + \imath \tau)) \right) - \varphi\left(\imath (\pi - \epsilon) \right) \nonumber \\
&= \frac{1}{4} \frac{\pi}{\epsilon} |z_l|^2 \, \left( \frac{1}{1 + \imath \tau} + \imath \tau - 1 \right) - \frac{1}{2} G^{''}(\epsilon) \epsilon^2 \tau^2 + K(\epsilon \tau), \label{p1p}
\end{align}
where
\begin{align} \label{ke}
K(\xi) = K(\mathcal{K}, \mathcal{A}; z, t; \xi) = G(\epsilon + \imath \xi) - G(\epsilon) - G'(\epsilon) \left( \imath \xi \right) - \frac{1}{2} G^{''}(\epsilon) \left( \imath \xi \right)^2.
\end{align}
We point out that $K(\xi)$ has been chosen so that
\begin{align*}
K(0) = K^{\prime}(0) = K^{\prime\prime}(0) = 0 \  \mbox{ and } \ K^{\prime\prime\prime}(0) = -\imath G^{\prime\prime\prime}(\epsilon).
\end{align*}

For the simplicity of notations, we write in the sequel
\begin{eqnarray} \label{ie1}
D_1 = \frac{1}{4} \frac{\pi}{\epsilon} |z_l|^2, \  D_2=\frac{1}{2} G^{\prime \prime}(\epsilon) \epsilon^2, \  J_* = K(\imath \epsilon) = G(0) - G(\epsilon) + G'(\epsilon) \epsilon - \frac{1}{2} G^{\prime\prime}(\epsilon) \epsilon^2.
\end{eqnarray}

This yields
\begin{eqnarray}
D_1 + D_2 = - \frac{1}{2} \epsilon^2 \Phi^{\prime\prime}(0),
\end{eqnarray}
and for $\epsilon > 0$ small enough
\begin{eqnarray} \label{d2j}
D_2 \sim |z|^2 \epsilon^2 > 0, \qquad J_* = O(|z|^2 \epsilon^3).
\end{eqnarray}

For $- \pi < \Im \omega < \pi$, we define
holomorphic functions $h_1$, $h_2$ and $h$ by
\begin{align}\label{defh}
h_1(\omega) = \frac{ \omega }{ \sinh{ \omega } }, \quad
h_2(\omega) = \omega \coth{ \omega }, \quad h(\omega) = h(\mathcal{K},\mathcal{A};\omega) = \prod_{j=1}^l h_1^{k_j} (a_j \omega).
\end{align}

Let $\omega = \xi + \imath \eta$ with $\xi, \eta \in \R $, then from a
direct computation, we get
\begin{eqnarray}\label{h1h2}
| h_1 ( \omega ) |^2 = \frac{ \xi^2 + \eta^2 }{ \sinh^2{ \xi } + \sin^2{ \eta } }, \qquad
\Re{h_2 ( \omega )} = \frac{ \xi \sinh{ \xi } \cosh{ \xi } + \eta \sin{ \eta } \cos{ \eta }}{ \sinh^2{ \xi } + \sin^2{ \eta } } \geq 0,
\end{eqnarray}
where $\Re \omega$ denotes the real part of a complex number $\omega$. \\

By \eqref{1}, \eqref{phi}, \eqref{de1} and the Cauchy's integral formula, we can write for $t \geq 0$ and $z_l \neq 0$
\begin{align} \label{spe}
p^{\H(\mathcal{K},\mathcal{A})}(z,t)
= \frac{1}{2(4\pi)^{n+1}}  e^{-\frac{d^2(z,t)}{4}} M, \quad \mbox{with} \quad M = \int_{\R} h(\lambda + \imath \theta) e^{\Phi(\lambda)} \, d\lambda.
\end{align}
Given $(z,t)$ such that $t \geq 0$ and $z_l \neq 0$, the above result is in fact a shift of the contour from
$\R$ to $\R + \imath \theta$, which can be derived from the exponential decay of the integrand in the horizonal strip bounded by $\R$ and $\R + \imath \theta$.\\

\subsection{Uniform asymptotic behaviour at infinity of the heat kernel at time $1$ on non-isotropic Heisenberg groups $\H(\mathcal{K}, \mathcal{A})$}\label{mainre}

Recall that $p^{\H(\mathcal{K},\mathcal{A})}$ denotes the heat kernel at time $1$ on the non-isotropic Heisenberg groups $\H(\mathcal{K}, \mathcal{A})$.
With the notation $\theta$, $d(z, t)$, $\Phi$, $\epsilon$, $D_1$, $D_2$, $J_*$, $G$, $G_1$, $G_2$ and $G_3$ defined above,
we will distinguish the following three cases:\\
 1. $d(z, t) \longrightarrow + \infty$ with $ | \theta| \leq \theta_0 < \pi$; \\
 2. $d(z, t) \longrightarrow + \infty$ with $0 < \epsilon = \pi - |\theta| \ll 1$ and $D_1 + D_2 \gg 1$;\\
 3. $d(z, t) \longrightarrow + \infty$ with $0 \leq \epsilon \ll 1$ and $D_1 + D_2$ bounded.  \\
 They are treated by next three theorems, respectively. Furthermore, we will find that these results are uniform. By the continuity of $p^{\H(\mathcal{K},\mathcal{A})}$,
 without loss of generality, we may assume that $|z| \neq 0$ and even $|z_l| \neq 0$ if we need.

\begin{theorem} \label{t1}
Suppose that $\frac{\pi}{2} \leq \theta_0 < \pi$. There exists a constant $c(\theta_0) = c(\mathcal{K}, \mathcal{A}; \theta_0) > 0$ such that:
\begin{align} \label{t1as}
p^{\H(\mathcal{K},\mathcal{A})}(z,t)
= \frac{e^{-\frac{d^2(z,t)}{4}}}{4 \, (4 \pi)^{n + \frac{1}{2}}} \left( -\frac{\Phi^{\prime\prime}(0)}{2} \right)^{-\frac{1}{2}} \prod_{j=1}^l \left( \frac{a_j \theta}{\sin{(a_j \theta)}} \right)^{k_j}
 \left( 1 + O_{\theta_0}(d^{-2}(z, t)) \right),
\end{align}
for all $d(z, t) \geq c(\theta_0)$ satisfying $|\theta| \leq \theta_0$.
\end{theorem}

\begin{theorem} \label{t2}
There exists a constant $0 < \epsilon_0 \ll 1$ such that
\begin{align}
p^{\H(\mathcal{K},\mathcal{A})}(z,t)
&= \frac{e^{-\frac{d^2(z,t)}{4}}}{4 \, (4 \pi)^{n + \frac{1}{2}}} (D_1 + D_2)^{-\frac{1}{2}} \, \epsilon
\prod_{j=1}^{l} \left( \frac{a_j \theta}{\sin{(a_j \theta)}} \right)^{k_j} \,   \left( 1 + O((D_1 + D_2)^{-1}) \right)  \label{t2as} \\
 &= \frac{e^{-\frac{d^2(z,t)}{4}}}{4 \, (4 \pi)^{n + \frac{1}{2}}}
\frac{\epsilon \left( \frac{\pi}{\epsilon} \right)^{k_l}}{\sqrt{D_1 + D_2}}
\prod_{j=1}^{l-1} \left( \frac{a_j \pi}{\sin{(a_j \pi)}} \right)^{k_j} \,   (1+O(\epsilon)) \left( 1 + O((D_1 + D_2)^{-1}) \right), \nonumber
\end{align}
for all $(z, t)$ satisfying
$0 < \epsilon = \pi - |\theta| \leq \epsilon_0$ and $D_1 + D_2 = - \frac{\Phi^{\prime\prime}(0)}{2} \epsilon^2 \ge \epsilon_0^{-3}$.
\end{theorem}

\begin{remark}
The order of the
remainders $O_{\theta_0}(d^{-2}(z, t))$ and $O((D_1 + D_2)^{-1})$, in Theorems \ref{t1} and \ref{t2} respectively, can be explained by the method of stationary phase (see for example \cite[Ch. VIII 1.3.4]{S93}).
\end{remark}

\begin{theorem} \label{t3}
Let $\epsilon_0$ be the constant in Theorem \ref{t2}. For $\gamma_0 \geq 1$ fixed, we have
\begin{eqnarray} \label{AT3}
p^{\H(\mathcal{K},\mathcal{A})}(z,t)
= \frac{e^{-\frac{d^2(z,t)}{4}}}{2 \, (4\pi)^{n+1}}
 e^{-D_1+J_*}  S_{k_l} \, \epsilon \left( \frac{\pi}{\epsilon} \right)^{k_l}
\prod_{j=1}^{l-1} \left( \frac{a_j \pi}{\sin{(a_j \pi)}} \right)^{k_j} \, \left( 1 + O_{\gamma_0} \left(d^{-\frac{1}{2}}(z, t)\right) \right),
\quad
\end{eqnarray}
as $d(z, t) \to +\infty$ satisfying $0 \leq \epsilon \leq \epsilon_0$ and $D_1 + D_2 \leq \gamma_0$,
with
\begin{align} \label{SKL}
 S_{k_l} = \frac{\sqrt{\pi}}{\sqrt{D_2}} \int_0^{+\infty} e^{-s}
\left( \frac{s}{D_1} \right)^{\frac{1}{2} (k_l - 1)} I_{k_l-1}(2 \sqrt{D_1 s}) e^{- \frac{(s - D_1)^2}{4D_2}} \, ds,
\end{align}
where the modified Bessel function $I_{\nu}$ ($\nu > -\frac{1}{2}$) is defined as
\begin{eqnarray*}
I_{\nu}(u)  =  \frac{1}{\sqrt{\pi} \Gamma(\nu + \frac{1}{2})} \Big(
\frac{u}{2} \Big)^{\nu} \int_{-1}^1 (1 - h^2)^{\nu - \frac{1}{2}}
e^{-u h} \, dh, \quad u \in \C.
\end{eqnarray*}
\end{theorem}

\begin{remark}
The remainder $O_{\gamma_0}(d^{-\frac{1}{2}}(z, t))$ in Theorem \ref{t3} should not be optimal and we should ask whether it can be improved to $O_{\gamma_0}(d^{-1}(z, t))$.
\end{remark}

\begin{remark}
When $l = 1$, that is, in the isotropic case, for $0 < \epsilon \ll 1$, a better asymptotic bound, which is given by a different formula, has been obtained in \cite[Lemme 2]{Li07}. A detailed explanation on the coincidence of their leading terms can be found in \cite{Li12}.
\end{remark}

Here we use the formula \eqref{AT3} in order to  match Theorem \ref{t3} with Theorem \ref{t2}. Someone may ask that the term $\epsilon^{1 - k_l} S_{k_l}$ in \eqref{AT3} could blow up for $k_l \geq 2$ as $\epsilon \longrightarrow 0^+$, and that \eqref{AT3} might not be uniform. In fact, by using the definition of $D_1$, $D_2$ and $\epsilon$ (see \eqref{ie1} and \eqref{epsilon}
above), the change of variables $s = \epsilon \rho$ shows that
\begin{align*}
\epsilon^{1 - k_l} S_{k_l} &= \sqrt{\pi} \left( \frac{1}{2} G^{''}(\epsilon) \right)^{-\frac{1}{2}} \\
&\times
\int_0^{+\infty} e^{-\epsilon \rho}
\left( \frac{4 \rho}{\pi |z_l|^2} \right)^{\frac{1}{2}(k_l-1)} I_{k_l-1}(\sqrt{\pi} |z_l| \sqrt{\rho}) \exp{\left(-\frac{(\rho - \frac{\pi |z_l|^2}{4 \epsilon^2})^2}
{2 G^{''}(\epsilon)}\right)} \, d\rho.
\end{align*}
Here if $|z| \neq 0$ and $\epsilon = 0$ (it implies $|z_l| = 0$), we have $G^{''}(0) = G_3^{''}(0)$ (cf. \eqref{Ge2}
above) and the following interpretations:
\begin{align*}
\lim_{|z_l| \longrightarrow 0} \left( \frac{4 \rho}{\pi |z_l|^2} \right)^{\frac{1}{2}(k_l-1)} I_{k_l-1}(\sqrt{\pi} |z_l| \sqrt{\rho}) = \frac{\rho^{k_l - 1}}{\Gamma(k_l)}, \quad \forall \rho > 0, \\
\lim_{|z_l| \longrightarrow 0} \frac{\pi |z_l|^2}{4 \epsilon^2} = \frac{1}{4} \left( |t| - \sum_{j=1}^{l-1} a_j \mu(a_j \pi) |z_j|^2 \right).
\end{align*}

In particular, we get the following corollary.

\begin{corollary} \label{cor1}
Suppose that $|z| \neq 0$ and $|\theta| = \pi$. Then $|z_l| = 0$ and we have
\begin{align}  \label{ATsT}
p^{\H(\mathcal{K},\mathcal{A})}(z,t)
&= \frac{\pi^{k_l - n - \frac{1}{2}}}{2^{2n + 3} \Gamma(k_l)}
 \left( \frac{1}{2} G^{''}_3(0) \right)^{-\frac{1}{2}} e^{-\frac{d^2(z,t)}{4}}
\prod_{j=1}^{l-1} \left( \frac{a_j \pi}{\sin{(a_j \pi)}} \right)^{k_j} \, \left( 1 + O\left(d^{-\frac{1}{2}}(z, t)\right) \right) \nonumber \\
&\times \int_0^{+\infty} \rho^{k_l - 1}
 \exp{\left(-\frac{1}{2 G^{''}_3(0)} \left( \rho - \frac{1}{4} ( |t| - \sum_{j=1}^{l-1} a_j \mu(a_j \pi) |z_j|^2 ) \right)^2
 \right)} \, d\rho.
\end{align}
\end{corollary}

The idea of the proof for Theorems \ref{t1} and \ref{t2} is standard, via the method of stationary phase. The main difficulty arises in the proof of Theorem \ref{t3}. Analytically, this is because $\pm \imath \pi$ are simultaneously poles of the phase and of the amplitude in the expression of the heat kernel by means of oscillatory integral (cf. \eqref{1} above), and the method of stationary phase is no longer valid. From a geometric point of view, this could be explained by the fact that the cut locus of $o$, excluding itself, is located in the case where $\epsilon = 0$. In addition, the phase (corresponding to the cut locus of $o$, or group law, from a geometric or algebraic point of view, respectively) in the non-isotropic case is much more complicated than that in the isotropic case, and the method in \cite{Li07} is not enough. To overcome this difficulty, we will use a more powerful method. More precisely we essentially follow the same strategy in \cite{Li12}:

First, we remark that a straightforward modification of \cite{Li12} suffices for $k_l \geq 2$ since the corresponding bad part of the phase for the case of Grushin operators in \cite{Li12} is more complicated than ours. Next, in the other case $k_l = 1$, corresponding to a case left open (in fact the case $n = 2$ there) in \cite{Li12}, the amplitude will play an important role as well, which makes the estimate of remainders subtler. Fortunately, it can be handled by ameliorating the proof of \cite[pp. 809-822]{Li12} (more specifically, by adding
a step of integration by parts in some appropriate form motivated by \cite[Cas 2, p. 816]{Li12}, so it should be reasonable to distinguish the case $k_l = 1$ with other cases $k_l \geq 2$).

\medskip

The idea and method of \cite{Li12} can be adapted to other situations. For example, we will provide an intrinsic method in a forthcoming work to prove the uniform asymptotic bounds at infinity for heat kernel in the setting of H-type groups, which have been obtained in \cite{Li10} by a very natural method motivated by \cite{Li07}. However, the method in \cite{Li10} is no longer valid for other stratified group $G$ satisfying $\mathrm{dim} \mathcal{V}_2 \geq 2$, where $\mathcal{V}_2$ denotes the second slice in the stratification of its Lie algebra.

\medskip

Moreover, we point out that we will use the explicit formulas for the Carnot-Carath\'eodory distance in the proof below. However, the readers, who are familiar with asymptotic methods, will find that it is not necessary. It is worth observing that the Carnot-Carath\'eodory distance can be obtained in some more general step-$2$ stratified groups, for example, Lie groups associated to quadratic CR manifolds in the sense of \cite{NRS01}.

\subsection{Precise estimates for the heat kernel and its derivatives on non-isotropic Heisenberg groups $\H(\mathcal{K}, \mathcal{A})$}
\label{precise}

A slight modification (substituting $\nu$ for $\frac{n}{2} - m \geq 1$, $r$ for $D_1$ and $b$ for $D_2$ respectively) of \cite[(4.11)]{Li12} and  \cite[Proposition 4.1, p. 817]{Li12} gives the following proposition, which plays a crucial role in \cite{Li12} and also in this article.

\begin{proposition}\label{prop1}
Let $\nu \ge 1$. For any $r \ge 0$ and $b > 0$,
\begin{align}\label{Plancherel}
\int_{\R} \frac{1}{(1+\imath \tau)^{\nu}} e^{-b \tau^2 + \imath r \tau + \frac{r}{1 + \imath \tau}} d\tau
= \frac{\sqrt{\pi}}{\sqrt{b}} \int_0^{\infty}
e^{-s} \left( \frac{s}{r} \right)^{\frac{\nu - 1}{2}} I_{\nu - 1}(2 \sqrt{r s})
 e^{- \frac{(s - r)^2}{4b}} \, ds.
\end{align}
Furthermore, for $\gamma_0 > 0$, we have
\begin{align} \label{KE1}
\frac{\sqrt{\pi}}{\sqrt{b}} \int_0^{\infty}
e^{-s} \left( \frac{s}{r} \right)^{\frac{\nu - 1}{2}} I_{\nu - 1}(2 \sqrt{r s})
 e^{- \frac{(s - r)^2}{4b}} \, ds \sim_{\gamma_0, \nu} \left(r + \sqrt{b}\right)^{\nu - 1},
\end{align}
for all $0 \le r \le \gamma_0$ and $0 < b \le \gamma_0$.
\end{proposition}

Using the scaling property \eqref{sp}, the definition of $S_{k_l}$ (cf. \eqref{SKL}) and \eqref{KE1} (with $\nu = k_l$, $r = D_1$ and $b = D_2$), we can conclude the following results.

\begin{corollary}
We have
\begin{eqnarray} \label{pehk}
p^{\H(\mathcal{K},\mathcal{A})}_h(z,t)
\sim h^{-n - 1} \left( 1 + \frac{|z|^2}{h} \epsilon_*^2 + \frac{|z_l|^2}{h} \epsilon_*^{-1} \right)^{- \frac{1}{2}}
 \left( \frac{h + \sqrt{h} |z| + |z_l|^2  \epsilon_*^{-2}}{h + \sqrt{h} |z| \epsilon_* + |z_l|^2 \epsilon_*^{-1}} \right)^{k_l - 1} e^{-\frac{d^2(z, t)}{4 h}}, \quad
\end{eqnarray}
for all $h > 0$ and $(z, t) = ((z_1, \cdots, z_l), t) \in \H(\mathcal{K}, \mathcal{A})$ satisfying $\epsilon_* = \frac{\sin{\theta}}{\theta} > 0$.
\end{corollary}

\begin{remark}
For the special case $\H((1, 1), (a_1, 1))$ (with $0 < a_1 < 1$), namely the five dimensional non-isotropic Heisenberg groups, \eqref{pehk} can be obtained by applying \eqref{HEB1} and the following counterpart of \eqref{convolution2}:
\begin{align}\label{convolutionh5}
p_h^{\H((1, 1), (a_1, 1))}(z, t) = \frac{1}{a_1} \left(p_h^{\H(1, 1)}\left(z_1, \frac{\cdot}{a_1}\right) \ast p_h^{\H(1, 1)}(z_2, \cdot)\right)(t).
\end{align}
However, the proof needs extremely tedious calculations and it seems not very practical for our purpose.
\end{remark}

\begin{remark}
By \eqref{dd2} above, we get a uniform, but more complicated, estimate as follows:
\begin{align*}
p^{\H(\mathcal{K},\mathcal{A})}(z,t)
& \ \sim  e^{-\frac{d^2(z, t)}{4}} \left[ 1 + |z|^2 \epsilon_*^2 + \epsilon_* \left( d^2(z, t) - \sum_{j = 1}^{l - 1} |z_j|^2 \left( \frac{a_j \theta}{\sin{(a_j \theta)}} \right)^2 \right) \right]^{- \frac{1}{2}} \\
&\times \left[ \frac{1 + |z| + \left( d^2(z, t) - \sum\limits_{j=1}^{l - 1}
|z_j|^2 \left( \frac{a_j \theta}{\sin{(a_j \theta)}} \right)^2 \right)}{1 + |z| \epsilon_* + \epsilon_* \left( d^2(z, t) - \sum\limits_{j=1}^{l - 1}
|z_j|^2 \left( \frac{a_j \theta}{\sin{(a_j \theta)}} \right)^2 \right)} \right]^{k_l - 1}, \quad \forall (z, t) \in \H(\mathcal{K}, \mathcal{A}).
\end{align*}
In the case where $l = 1$, that is, in the isotropic Heisenberg group case, it becomes (see also \eqref{IHE})
\begin{eqnarray*}
p^{\H(n , 1)}(z,t) \sim \left( 1 + |z| d(z, t) \right)^{- \frac{1}{2}}
 \left( \frac{1 + d^2(z, t)}{1 + |z| d(z, t)} \right)^{n - 1} e^{-\frac{d^2(z, t)}{4}}, \qquad \forall (z, t) \in \H(n, 1),
\end{eqnarray*}
which has been obtained in \cite[Theorem 1]{Li07}. However, any formula similar to the last one, depending only on $(|z|, d(z, t))$, is invalid in the nonisotropic case.
\end{remark}

Furthermore, we give sharp estimates for derivatives of $p^{\H(\mathcal{K},\mathcal{A})}$ (so of $p^{\H(\mathcal{K},\mathcal{A})}_h$ as well) as follows:

\begin{theorem} \label{t4}
We have
$| \nabla_{\H(\mathcal{K},\mathcal{A})} \ln{p^{\H(\mathcal{K},\mathcal{A})}(g)}
| \lesssim d(g)$ for all $g \in \H(\mathcal{K}, \mathcal{A})$. In
general, let $m \in \N^+$, $\beta \in \N^m$ and $\U_{\eta} \in \{\X_{i, j}, \Y_{i, j} \}_{1 \leq i \leq l, 1 \leq j \leq k_i}$ with $1 \leq \eta \leq m$. Then there
exists a constant $C_* = C_*(\beta) > 0$ such that
\begin{align*}
| \U_1^{\beta_1} \cdots \U_m^{\beta_m} p^{\H(\mathcal{K},\mathcal{A})}(g) | \leq C_* (1 + d(g))^{|\beta|} p^{\H(\mathcal{K},\mathcal{A})}(g), \quad \forall g \in \H(\mathcal{K}, \mathcal{A}), \
\mbox{where }|\beta| = \sum_{\eta = 1}^m \beta_{\eta}.
\end{align*}
\end{theorem}

\begin{remark}
The estimate $| \nabla_{\H(\mathcal{K},\mathcal{A})} \ln{p^{\H(\mathcal{K},\mathcal{A})}(g)} | \lesssim d(g)$ has been obtained by B. Qian in \cite{Q10} by means of Harnack inequality and Bakry-Emery criterion. See \cite{BB09} for more details about the method.
\end{remark}

Notice that we can also obtain the
asymptotic behaviour at infinity of $\U_1^{\beta_1} \cdots \U_m^{\beta_m} p^{\H(\mathcal{K},\mathcal{A})}$.
They are left to the interested reader.

In a future work, the second author will establish the gradient bound of the heat semigroup on $\H(\mathcal{K}, \mathcal{A})$.

\subsection{Small-time asymptotic behaviour on $\H(\mathcal{K}, \mathcal{A})$ of the heat kernel}

Since the result of the case $(z, t) = o$ is trivial, we may assume in the sequel $(z, t) \neq o$. To study the small-time asymptotic estimates, Theorem \ref{t1} and \eqref{ATsT} are sufficient. In fact, from the scaling invariance property, Theorem \ref{t1} implies the following corollary.

\begin{corollary} \label{c2}
Suppose that $|t| = \sum\limits_{j=1}^{l}
a_j \mu(a_j \theta) |z_j|^2$ has a unique solution $0 \leq \theta < \pi$, then we have
\begin{eqnarray*}
p^{\H(\mathcal{K},\mathcal{A})}_h(z,t)
= \frac{1}{4 (4\pi h)^{n + \frac{1}{2}}} e^{-\frac{d^2(z,t)}{4h}} \left( -\frac{\Phi^{\prime \prime}(0)}{2} \right)^{-\frac{1}{2}} \prod_{j=1}^l \left( \frac{a_j \theta}{\sin{(a_j \theta)}} \right)^{k_j}
 \left( 1 + O(h) \right), \quad  \mbox{as $h \to 0^+$.}
 \end{eqnarray*}
\end{corollary}

By \eqref{ATsT}, we get the following corollary immediately.

\begin{corollary} \label{c3}
Assume that $z_l = 0$ and $|t| = \sum\limits_{j=1}^{l - 1}
a_j \mu(a_j \pi) |z_j|^2 > 0$. Then
\begin{eqnarray*}
p^{\H(\mathcal{K},\mathcal{A})}_h(z,t)
= \frac{2^{\frac{k_l - 3}{2}} \pi^{k_l + \frac{1}{2}}}{(4\pi)^{n+1}} \frac{\Gamma(\frac{k_l}{2})}{\Gamma(k_l)} \left( G_3^{\prime \prime}(0) \right)^{\frac{k_l - 1}{2}}  e^{-\frac{d^2(z,t)}{4h}}  h^{- n - \frac{k_l + 1}{2}}
 \prod_{j = 1}^{l-1} \left( \frac{a_j \pi}{\sin{(a_j \pi)}} \right)^{k_j}
 (1+o(1)),
\end{eqnarray*}
as $h \to 0^+$.
\end{corollary}

For $|z| \neq 0$, using a change of variables $$\rho = 4^{-1} \left( |t| - \sum_{j=1}^{l-1} a_j \mu(a_j \pi) |z_j|^2  \right) + \sqrt{2 G_3^{''}(0)} \lambda,$$
we obtain from \eqref{ATsT} the following corollary as well.

\begin{corollary} \label{c4}
If $z_l = 0$ and $|t| > \sum\limits_{j=1}^{l - 1} a_j \mu(a_j \pi) |z_j|^2$. Then
\begin{align*}
p^{\H(\mathcal{K},\mathcal{A})}_h(z,t)
&= \frac{4^{1 - k_l} \pi^{k_l + 1}}{(4 \pi)^{n+1} \Gamma(k_l)}
\left( |t| - \sum_{j=1}^{l-1} a_j \mu(a_j \pi) |z_j|^2  \right)^{k_l - 1} e^{-\frac{d^2(z,t)}{4h}} h^{ - n - k_l} \\
&\times \prod_{j=1}^{l-1} \left( \frac{a_j \pi}{\sin{(a_j \pi)}} \right)^{k_j} (1+o(1)), \end{align*}
as $h \to 0^+$.
\end{corollary}

\begin{remark}
Note that we have assumed that $|z| \neq 0$ in Corollary \ref{cor1}, but in the case of $z = o$ and $|t| > 0$, it suffices to use an argument of limit.
\end{remark}

\begin{remark}
Corollary \ref{c2} and Corollary \ref{c4} have been obtained by Beals, Gaveau and Greiner, namely \cite[Theorem 5.62 and Theorem 5.75]{BGG00} (with $u = 2 h$), respectively, but
Corollary \ref{c3} is new, which  completes their results. In particular, some new information, about the power of $h$ and the coefficient of the leading term, can be read from our result.
\end{remark}

\begin{remark}
The main results of \cite{BBN16}, based on Laplace's method, can be deduced from our Corollaries. In \cite{BBN16}, the authors considered only the very special situation of five dimensional non-isotropic Heisenberg groups and only $\sigma(g)$, namely, the power of $h$ in the term $\Theta_h(g) = C(g) h^{-\sigma(g)}$ of \eqref{STA} and the essential information can not be revealed.
\end{remark}

\medskip

As we pointed out above, this paper is closely related to \cite{Li12} and that there exist two cases left open in \cite{Li12}. So we will  fill in the blanks.

\subsection{The heat kernel for Grushin operators $\Delta_G = \sum_{i = 1}^{n} \left( U_i^2 + V_i^2 \right)$ (with $U_i = \frac{\partial}{\partial x_i}$ and $V_i = x_i \frac{\partial}{\partial t}$) on $\R^n \times \R$ w.r.t. Lebesgue measure $dx \, dt$} \label{Grushin}

Recall that the Grushin operators are prototype of nil-manifolds, but they do not have group structure. However, there exists dilatation structure as well. In \cite{Li12}, the first author has studied uniform asymptotic behaviour at infinity of the heat kernel at time $1$ as well as its precise estimates and its small-time asymptotic behaviour, although an additional condition, $n \geq 3$, is needed for main results therein, e.g. Theorems 1.2, 1.3 and Corollary 1.4 in Section 1.2, results in Sections 1.3, 1.4 and 1.5 and Proposition 4.2.

Now, we are in the place to remove this assumption and finish the remaining case in \cite{Li12}
by the trick of integration by parts used in the estimates of $M_2$ and $M_2^*$ in Subsections \ref{S52} and \ref{S52*} below (see Section \ref{st7} for more details). As a result, we have the following theorem.

\begin{theorem}\label{g}
Theorems 1.2 and 1.3 of \cite{Li12} are in fact valid for all $n \ge 1$.
\end{theorem}

\begin{remark}
By Theorem \ref{g}, Corollary 1.4 in Section 1.2, results in Sections 1.3, 1.4 and 1.5 and Proposition 4.2 of \cite{Li12} are still valid for $n = 1 $ and $n = 2$.
The remainders in Theorem 1.1 and Theorem 1.2 of \cite{Li12} can be improved to $O(d^{-2})$ and $O(\epsilon) + O((D_1 + D_2)^{-1})$ respectively
for all $n \ge 1$ if we use third order Taylor expansion as in the proof of Theorem \ref{t1} and Theorem \ref{t2} in this paper.
\end{remark}

Notice that in \cite{BBN12}, Barilari et al. have studied the small-time heat kernel asymptotics on Grushin plane, that is in the case where $n = 1$. Recall that their main result (valid in other situations, but under some technical assumptions) is proved by using only the standard Laplace's method, which is a special case of the method of stationary phase. However, as we mentioned above, even the method of stationary phase is not powerful enough for our purpose:
uniform asymptotic behaviour at infinity of the heat kernel.

\subsection{Further results}\label{futher}

It is worthwhile to point out that non-isotropic Heisenberg groups are isomorphic to the isotropic ones. In other words, instead of working on
$(\H(\mathcal{K},\mathcal{A}), \Delta_{\H(\mathcal{K},\mathcal{A})})$, we  can equivalently work on $\H(n,1)$, but with the sub-Laplacian
\begin{align*}
\sum_{i = 1}^l a_i \sum_{j = 1}^{k_i} (\hat{\X}_{i, j}^2 + \hat{\Y}_{i, j}^2),
\end{align*}
where
\begin{eqnarray*}
\hat{\X}_{i, j} = \frac{\partial}{\partial u_{i, j}} + 2  v_{i, j} \frac{\partial}{\partial t},  \quad \hat{\Y}_{i, j} = \frac{\partial}{\partial
v_{i, j}} - 2  u_{i, j} \frac{\partial}{\partial t}, \qquad 1 \leq i \leq l, 1 \leq j \leq k_i.
\end{eqnarray*}

More generally, let $G$ be a stratified Lie group with Lie algebra $\mathfrak{g}$ satisfying $\mathrm{dim} [\mathfrak{g},\mathfrak{g}] = 1$,
and $L$ a (not necessarily the canonical)  sub-Laplacian. By \cite[Ch.\,3 and Ch.\,16]{BLU07}, we can show that
there exist $q \in \{0, 1, \cdots\}$ and $(\mathcal{K}, \mathcal{A})$ such that $(G, L)$ is isomorphic to $\R^q \times \H(\mathcal{K}, \mathcal{A})$ with the canonical sub-Laplacian $\Delta_{\R^q \times \H(\mathcal{K}, \mathcal{A})}
= \Delta_{\R^q} + \Delta_{\H(\mathcal{K}, \mathcal{A})}$, and the sense for the case where $q =0$ is clear. By the fact that
\begin{align*}
e^{h \Delta_{\R^q \times \H(\mathcal{K}, \mathcal{A})}} = e^{h \Delta_{\R^q}} e^{h \Delta_{ \H(\mathcal{K}, \mathcal{A})}},
\end{align*}
the asymptotic estimates at infinity for the heat kernel associated to $(G, L)$ is now well-known.

\subsection{Structure of the paper}

The paper is organized as follows.  Theorem \ref{TC} will be proved in Section \ref{st1}. As we mentioned above that the proof of Theorems \ref{t1}-\ref{t3} in the case where $k_l \geq 2$ follows from a modification of the proof in \cite{Li12}. Then we only consider $k_l = 1$ here, which will be treated in Sections \ref{st1n}-\ref{st3}, respectively. In Section \ref{st4}, we give the proof of precise bounds for the heat kernel and its derivatives. An interpretation for the Grushin heat kernel in 2D and 3D will be presented in Section \ref{st7}.

\medskip

\renewcommand{\theequation}{\thesection.\arabic{equation}}
\section{Proof of Theorem 1} \label{st1}
\setcounter{equation}{0}

In this section we always assume
$l = 1$ and $n \geq 2$. By the scaling property, cf. \eqref{sp}, it suffices to show \eqref{IHE} for $h = 1$.

In such a case, according to \eqref{dd1} and \eqref{dd2}, we have for $|z| \neq 0$,
\begin{align*}
d^2(z, t) = |z|^2 \left( \frac{ \theta }{ \sin{ \theta } } \right)^2 \quad \mbox{with} \quad t = |z|^2 \mu(\theta) = |z|^2 \frac{2 \theta - \sin{2 \theta}}{2 \sin^2{\theta}}.
\end{align*}
Then
\begin{align*}
\frac{\partial \theta}{\partial t} = \frac{1}{|z|^2 \mu'(\theta)},
\end{align*}
and combining it with the fact that
\begin{align*}
\mu^{\prime} ( \theta ) = \frac{ 2( \sin{ \theta }-\theta \cos{ \theta } ) }{\sin^3{ \theta }}, \quad  \quad \frac{ d }{ d\theta } \left( \frac{ \theta }{ \sin{ \theta } } \right)^2 = \theta \mu'(\theta),
\end{align*}
we get the following lemma.

\begin{lemma}\label{l2}
For $ z \ne 0 $, $ \frac{ d }{ dt }d^2(z,t) = \theta $. As a consequence, $d(z, \cdot)$ is an even function and strictly increasing on the interval $(0, +\infty)$.
\end{lemma}

By an abuse of notation, let $d$ simultaneously denote the Carnot-Carath\'eodory distance on $\H(n, 1)$ and on the three dimensional Heisenberg group $\H(1, 1)$ in this section. Then, by \eqref{dd2}, we can write
\begin{align*}
d(z, t) = d(|z|, |t|), \qquad \forall (z, t) \in \H(n, 1).
\end{align*}

Furthermore, we will use the following simple fact:

\begin{lemma}
For $|z| \ne 0$, we have
\begin{eqnarray}\label{eps}
\epsilon = \pi - |\theta|
= \sqrt{ \pi } \sqrt{ \frac{ |z|^2 }{ |t| }}
\left( 1 + O \left( \sqrt{ \frac{ |z|^2 }{ |t| } } \right) \right)\qquad
 \mbox{as $\frac{ |t| }{ |z|^2 } \to +\infty$.}
\end{eqnarray}
\end{lemma}

\noindent \textbf{Proof.}
It suffices to use the fact that
\begin{eqnarray*}
\frac{ |t| }{ |z|^2 } = \mu ( |\theta| ) = \frac{ ( |\theta| - \sin{|\theta|}\cos{|\theta|} )}{ \pi } \frac{ (\pi - |\theta| )^2}{\sin^2{ |\theta| }} \frac{ \pi }{ (\pi - |\theta| )^2 }
=  \frac{ \pi }{ \epsilon^2 } (1 + O( \epsilon ))
\end{eqnarray*}
as $\frac{ |t| }{ |z|^2 } \to +\infty$ (it implies $\epsilon \to 0$) by \eqref{dd1}.
~ \hspace*{20pt} ~ \hfill $\Box$

The following result can be found in \cite[Th\'eor\`eme 2, p. 115]{G77}, or obtained by a direct calculation:

\begin{lemma}
We have
{\em \begin{align*}
\widehat{h_1}(s) := \int_{\R} \frac{\lambda}{\sinh{\lambda}} e^{-\imath \lambda s} \, d\lambda = \frac{\pi^2}{1 + \cosh{(\pi s)}} \sim e^{-\pi |s|}, \qquad \forall s \in \R.
\end{align*}}
where $\widehat{h_1}$ denotes the  Fourier transform of $h_1$ (cf. \eqref{defh}).
\end{lemma}

Using the following elementary property of Fourier transform and convolution:
\begin{align}\label{Fourier}
\widehat{f_1 \cdot f_2} = \frac{1}{2 \pi} \widehat{f_1} \ast \widehat{f_2}, \qquad \mbox{for suitable $f_1$ and $f_2$,}
\end{align}
an easy induction argument gives the following lemma.

\begin{lemma} \label{Lm4}
For $n \geq 2$, we have
\begin{align*}
\widehat{h_1^{n - 1}}(s) = \int_{\R} \left( \frac{\lambda}{\sinh{\lambda}} \right)^{n - 1} e^{-\imath \lambda s} \, d\lambda \sim (1+ \pi |s|)^{n - 2} e^{-\pi |s|}, \qquad \forall s \in \R.
\end{align*}
\end{lemma}

In the sequel of this section, we use $p(z,t)$ to denote $p^{\H(n, 1)}(z, t)$ for simplicity. Recall that without loss of generality, we may assume that $t \geq 0$, $|z| \ne 0$ and $d^2(z, t)$ is large enough.

From \eqref{convolution}, we can write
\begin{align*}
p(z,t) & =  \left( p^{\H(1, 1)}(|z|, \cdot) \ast p^{\H(n - 1, 1)}(o, \cdot) \right)(t) \sim \int_{\R}  \widehat{h_1^{n - 1}}\left(\frac{s}{4}\right) p^{\H(1, 1)}(|z|, t - s) ds \\
&\sim \int_{\R} ( 1 + \pi |s| )^{n-2} e^{- \frac{ \pi |s| }{ 4 }} e^{-\frac{ d^2(z , t-s )}{ 4 }}
\frac{1}{( 1 + |z| d(z, t-s))^{\frac{1}{2}}}  \, ds = \int_0^{+\infty} + \int_{-\infty}^0,
\end{align*}
where the last ``$\sim$'' follows from \eqref{HEB1} and Lemma \ref{Lm4}.
Now Lemma \ref{l2} implies that
\begin{align*}
d(z, t - s) \leq d(z, t + s), \qquad \forall t, s \geq 0,
\end{align*}
which yields that
\begin{align*}
p(z, t)  \sim \int_0^{+\infty} ( 1 + \pi s )^{n-2} e^{- \frac{ \pi s }{ 4 }} e^{-\frac{ d^2(z , t-s )}{ 4 }}
\frac{1}{( 1 + |z| d(z, t-s))^{\frac{1}{2}}}  \, ds=e^{-\frac{d^2(z,t)}{4}} Q,
\end{align*}
with
\begin{align}
Q =  \int_0^{+\infty} ( 1 + \pi s )^{n-2} e^{- \frac{ \pi s + d^2(z ,  t - s ) - d^2(z ,  t )}{ 4 }}
\frac{1}{( 1 + |z| d(z, t-s))^{\frac{1}{2}}}  \, ds.
\end{align}

It remains to show that
\begin{align} \label{eQ}
Q \sim  \left( 1 + |z| d(z, t) \right)^{-\frac{1}{2}} \left(   \frac{d^2(z, t)}{1 + |z| d(z, t)}  \right)^{n - 1}, \qquad \forall d(z, t) \gg 1.
\end{align}
We split the proof into the cases: (1) $\frac{t}{|z|^2} \geq \gamma_0$ with $|z| d(z,t) > 1$, (2) $\frac{t}{|z|^2} \geq \gamma_0$ with $|z|d(z,t) \leq 1$, and (3) $\frac{t}{|z|^2} \leq \gamma_0$, where the constant $\gamma_0$ is large enough ($\gamma_0 \ge 1000$ will be sufficient).

Case (1):  $\frac{t}{|z|^2} \geq \gamma_0$ with $|z| d(z,t) > 1$. We split the integral $Q$ into two parts as follows:
\begin{eqnarray*}
Q = \int_0^{ \frac{t}{2} } + \int_{ \frac{t}{2} }^{ +\infty } = Q_1 + Q_2.
\end{eqnarray*}
It follows from \eqref{l1} that $d(z, t - s) \sim \sqrt{t} \sim d(z, t)$ for all $0 < s \leq \frac{t}{2}$, then
\begin{align*}
Q_1 \sim \left( 1 + |z| d(z, t) \right)^{-\frac{1}{2}} \int_0^{ \frac{t}{2} }  ( 1 + \pi s )^{n-2} e^{- \frac{ \pi s + d^2(z ,  t - s ) - d^2(z ,  t )}{ 4 }} \, ds.
\end{align*}

Using the mean value theorem, Lemma \ref{l2} implies that for $0 < s \leq \frac{t}{2}$ we have
\begin{align*}
d^2(z, t) - d^2(z, t-s) = \theta_{\zeta_s} s,
\end{align*}
where $\theta_{\zeta_s} = \mu^{-1} \left( \frac{\zeta_s}{|z|^2} \right)$
with $\zeta_s \in (\frac{t}{2}, t)$. Under the assumption that $\frac{t}{|z|^2} \geq \gamma_0 \gg 1$, it deduces from the equation above and \eqref{eps} that

\begin{align*}
\frac{\pi s - d^2(z,t) + d^2(z,t-s)}{4} = \frac{s}{4}( \pi - \theta_{\zeta_s})
= \frac{s}{4} \sqrt{\pi} \sqrt{\frac{|z|^2}{\zeta_s}}(1 + o(1)).
\end{align*}
As a result, there exist two constants $C_1, C_2 > 0$ such that
\begin{align}\label{mean1}
C_2 \sqrt{\frac{|z|^2}{t}} s \leq
\frac{\pi s - d^2(z,t) + d^2(z,t-s)}{4}
\leq C_1 \sqrt{\frac{|z|^2}{t}} s.
\end{align}

Hence, we have
\begin{align*}
Q_1 &\lesssim \frac{ 1 }{( 1 + |z| d(z,t))^{\frac{1}{2}}} \int_0^{ \frac{t}{2} }
(1 + \pi s)^{n - 2} e^{ - C_2 \sqrt{ \frac{|z|^2}{t}} s} ds \\
&\sim  \frac{ 1 }{( 1 + |z| d(z, t))^{\frac{1}{2}}}
\left( \sqrt{\frac{t}{|z|^2} }\right)^{n - 1}
\sim \left( 1 + |z| d(z, t) \right)^{-\frac{1}{2}} \left(   \frac{d^2(z, t)}{1 + |z| d(z, t)}  \right)^{n - 1},
\end{align*}
where the second ``$\sim$'' follows from $d(z,t) \sim \sqrt{t}$ by \eqref{l1} and $|z| d(z,t) > 1$ and the first
``$\sim$'' follows from the change of variables $u = \sqrt{\frac{|z|^2}{t}}s$ in addition.

On the other hand, similarly, we have
\begin{align*}
Q_1 &\gtrsim \frac{ 1 }{( 1 + |z| d(z, t))^{\frac{1}{2}}} \int_0^{ \frac{t}{2} }
(1 + \pi s)^{n - 2} e^{ - C_1 \sqrt{ \frac{|z|^2}{t}} s} ds \\
&\sim  \frac{ 1 }{( 1 + |z| d(z, t))^{\frac{1}{2}}}
\left( \sqrt{\frac{t}{|z|^2} }\right)^{n - 1}
\sim \left( 1 + |z| d(z, t) \right)^{-\frac{1}{2}} \left(   \frac{d^2(z, t)}{1 + |z| d(z, t)}  \right)^{n - 1}.
\end{align*}

To estimate $Q_2$, we use the mean value theorem again to obtain
\begin{align*}
d^2(z, t) - d^2(z, t-s) = \theta_{\zeta_s} s \leq \theta_t s,  \qquad \forall s \geq \frac{t}{2},
\end{align*}
where $\theta_{\zeta_s} = \mu^{-1} \left( \frac{\zeta_s}{|z|^2} \right)$
with $\zeta_s \in (t - s, t)$ and $\theta_t = \mu^{-1} \left( \frac{t}{|z|^2} \right)$. As before, the above inequality and \eqref{eps} imply that there exists a constant $C_3 > 0$ such that
\begin{eqnarray}\label{mean2}
-\frac{\pi s - d^2(z, t) + d^2(z, t - s)}{4}
\leq - \frac{s}{4} (\pi - \theta_t)
\leq -C_3 \sqrt{\frac{|z|^2}{t}} s.
\end{eqnarray}

Then we have
\begin{align*}
Q_2  \le \int_{ \frac{t}{2} }^{ +\infty }
(1 + \pi s)^{n - 2} e^{ - C_3 \sqrt{ \frac{|z|^2}{t}} s} \, ds \sim \int_{ \frac{t}{2} }^{ +\infty } s^{n - 2} e^{ - C_3 \sqrt{ \frac{|z|^2}{t}} s} \, ds.
\end{align*}
By the change of variables $u = \sqrt{ \frac{|z|^2}{t}} s$ in the last integral, we yield
\begin{align*}
Q_2 &\lesssim \left( \sqrt{\frac{t}{|z|^2} }\right)^{n - 1}
\int_{ \frac{\sqrt{|z|^2t} }{2}}^{ +\infty }  u^{n - 2}  e^{ - C_3u } \, du \lesssim  \left( \sqrt{\frac{t}{|z|^2} }\right)^{n - 1}
e^{- \frac{C_3}{2} \frac{ \sqrt{|z|^2 t }}{2}}.
\end{align*}
Recall that $\sqrt{|z|^2 t } = |z| \sqrt{t} \sim |z| d(z, t) \geq 1$ which yields that $e^{- \frac{C_3}{2} \frac{ \sqrt{|z|^2 t }}{2}} \lesssim \left( 1 + |z| d(z, t) \right)^{-\frac{1}{2}}$. Hence we obtain the validity of \eqref{eQ} in the case $\frac{t}{|z|^2} \geq \gamma_0  \gg 1$ with $|z| d(z,t) > 1$.

\vskip4pt

Case (2):  $\frac{t}{|z|^2} \geq \gamma_0$ with $|z| d(z,t) \leq 1$.
Now we split the integral $Q$ into three parts as follows:
\begin{eqnarray*}
Q = \int_0^{ \frac{t}{2} } + \int_{ \frac{t}{2} }^{ 2t } + \int_{ 2t }^{ +\infty } = Q_1 + Q_{2,1} + Q_{2,2}.
\end{eqnarray*}

Since $\frac{t}{|z|^2} \geq \gamma_0 \gg 1$, we have $d^2(z, t) \sim t$ is large enough from \eqref{l1}. Using
\eqref{mean1}, \eqref{mean2} again and the fact that
\begin{align*}
\sqrt{\frac{|z|^2}{t}} s \leq 2 |z| \sqrt{t} \lesssim |z| d(z, t)  \le 1, \qquad \forall 0 \leq s \leq 2 t,
\end{align*}
from our assumption, we may easily obtain
\begin{align*}
Q_1 \sim \int_{0}^{ \frac{t}{2} } (1 + \pi s)^{n - 2} \, ds \sim t^{n - 1},
\end{align*}
and
\begin{eqnarray*}
Q_{2,1} \lesssim \int_{ \frac{t}{2} }^{ 2t } (1 + \pi s)^{n - 2} ds \sim t^{n - 1}.
\end{eqnarray*}

To estimate $Q_{2,2}$, observe that it deduces from Lemma \ref{l2} that
\begin{align*}
d(z,t) \leq d(z, t - s) = d(z, s- t), \qquad \forall s \geq 2 t,
\end{align*}
which implies that
\begin{eqnarray*}
Q_{2,2}  \le \int_{ 2t }^{ +\infty } (1 + \pi s)^{n - 2} e^{ -\frac{\pi}{4} s} \, ds = o\left(e^{-\frac{t}{4}}\right).
\end{eqnarray*}
This proves the required estimate \eqref{eQ} in the case $\frac{t}{|z|^2} \geq \gamma_0$ with $|z| d(z,t) \leq 1$
since under the assumption we have $|z| d(z,t) + 1 \sim 1$ and $d^2(z,t) \sim t$ from \eqref{l1}.

\vskip4pt

Case (3): $\frac{t}{|z|^2} \leq \gamma_0$. Note that by \eqref{l1} and Lemma \ref{l2},
\begin{align}\label{ppd}
d^2(z, t) \sim_{\gamma_0} |z|^2 \gg 1, \qquad d(z, t) \lesssim_{\gamma_0} d(z, 0) \leq d(z, s), \quad \forall s \in \R.
\end{align}
Using Lemma \ref{l2} and the mean value theorem as before, we get for any $s > 0$
\begin{align*}
d^2(z,t) - d^2(z, t - s) = \theta_{\zeta_s} s,
\end{align*}
where $\theta_{\zeta_s} = \mu^{-1}\left( \frac{\zeta_s}{|z|^2} \right) \in (-\pi, \theta_0]$ and $\theta_0=\mu^{-1}(\gamma_0)$ is a constant belonging to $(0,\pi)$. Then
\begin{eqnarray}\label{mean3}
- \frac{2 \pi}{4} s \leq - \frac{\pi s - d^2(z,t) + d^2(z,t-s)}{4}
\leq - \frac{\pi - \theta_0}{4} s, \quad \forall s > 0.
\end{eqnarray}

On one hand, by the second inequality of \eqref{mean3} and \eqref{ppd} we have
\begin{align*}
Q &= \int_0^{+\infty}
(1 + \pi s)^{n - 2} \frac{ e^{ -\frac{\pi s -  d^2(z,t) + d^2(z, t - s)}{4}} }{( 1 + |z| d(z, t-s))^{\frac{1}{2}}} \, ds\\
& \lesssim \frac{1}{|z|} \int_0^{+\infty} (1 + \pi s)^{n - 2} e^{- \frac{\pi - \theta_0}{4} s} ds \lesssim_{\gamma_0} \frac{1}{|z|}.
\end{align*}
On the other hand,
\begin{align*}
Q \gtrsim_{\gamma_0} \frac{1}{|z|} \int_0^1 (1 + \pi s)^{n - 2} e^{-\frac{\pi}{2} s} \, ds \sim \frac{1}{|z|}
\end{align*}
by the first inequality of \eqref{mean3} and \eqref{l1}. The two estimates above and \eqref{ppd} complete the proof in this case.
Combining all results above, we obtain \eqref{eQ}.
~ \hspace*{20pt} ~ \hfill $\Box$

\medskip

\renewcommand{\theequation}{\thesection.\arabic{equation}}
\section{Uniform asymptotic behaviour of $M$ (cf. \eqref{spe}) for $d(z, t) \longrightarrow +\infty$ with $0 \leq \theta \leq \theta_0 < \pi$ when $l \geq 2$ and $k_l = 1$}  \label{st1n}
\setcounter{equation}{0}

In this paper, we need the following lemma which is a counterpart of (4.5) in \cite{Li10}:

\vskip4pt

\begin{lemma}\label{l5}
$\Re\Phi(s)$ (cf. \eqref{phi}) is an even function independent of $t$, strictly decreasing on the interval
$[0,+\infty)$ and there exists a constant $A_0 = A_0(\mathcal{K},\mathcal{A}) > 0$ such that
\begin{eqnarray}\label{RePhi}
\Re\Phi(s) \leq -A_0 |z|^2 s^2, \qquad  \forall -\pi < \theta <\pi, 0< s \leq1.
\end{eqnarray}
\end{lemma}

\vskip4pt

 \noindent \textbf{Proof.}
From \eqref{defvarphi}, \eqref{phi} and \eqref{h1h2}, we have
\begin{align*}
\Re\Phi(s) &= -\frac{ t \theta }{ 4 }- \sum_{j=1}^{l} \frac{ |z_j|^2 }{ 4 }
\frac{ a_j s \sinh( a_j s) \cosh( a_j s ) + a_j \theta \sin( a_j \theta )\cos( a_j \theta )}
{\sinh^2( a_j s )+\sin^2( a_j \theta )} \\
&\quad +\frac{ t \theta }{ 4 } + \sum_{j=1}^{l} \frac{ |z_j|^2 }{ 4 } a_j \theta \cot( a_j \theta )\\
&=-\sum_{j=1}^{l} \frac{ |z_j|^2 }{ 4 } \frac{ \sinh^2( a_j s )}{ \sinh^2( a_j s ) + \sin^2( a_j \theta )}
( a_j s \coth( a_j s ) - a_j \theta \cot( a_j \theta )),
\end{align*}
and our lemma follows.
~ \hspace*{20pt} ~ \hfill $\Box$

\vskip4pt

In the rest of the paper, we use $p(z,t)$ and $\nabla$ to denote $p^{\H(\mathcal{K}, \mathcal{A})}(z, t)$ and $\nabla_{\H(\mathcal{K}, \mathcal{A})}$ respectively for simplicity.

Recall that we can assume that $t \geq 0$ and $\frac{\pi}{2} \leq \theta_0 < \pi$. Under our assumption,
\eqref{dep} and \eqref{dd2} imply that
\begin{align} \label{ep11}
-\Phi^{\prime\prime}(0) \sim_{\theta_0} |z|^2 \sim_{\theta_0} d^2(z, t) \gg 1.
\end{align}

Observe that, see \eqref{spe},
\begin{align*}
M =  \int_{\R} h(\lambda + \imath \theta) e^{\Phi(\lambda)} \, d\lambda
= \int_{-\delta}^{\delta} + \int_{|\lambda| > \delta } = T_1 + T_2,
\end{align*}
where the constant $0 < \delta = \delta(\mathcal{K}, \mathcal{A}; \theta_0; z) \ll 1$ will be chosen later. By the basic property of heat kernel, we have $p(z, t) > 0$ and it follows from \eqref{spe} that $M > 0$. Hence, we will control in the sequel the moduli of negligible terms but don't
care about its imaginary part.

Noting that the function $h(\omega) = \prod\limits_{j=1}^l h_1^{k_j}(a_j \omega)$ is holomorphic for $-\pi < \Im \omega < \pi$ (see the definition of $h$ and $h_1$ in \eqref{defh}) and by the fact that $h(\imath \theta) \geq 1$, we have
\begin{align} \label{1h1}
h(\lambda + \imath \theta) = h( \imath \theta ) \left( 1 + \frac{h^{\prime}( \imath \theta )}{h( \imath \theta )} \lambda  + O_{\theta_0}( \lambda^2 ) \right), \quad \forall - \frac{\pi - \theta_0}{4} \leq \lambda \leq \frac{\pi - \theta_0}{4}, \  0 \leq \theta \leq \theta_0.
\end{align}

Moreover, it deduces from the definition of  $\Phi$ and Taylor's formula (cf. \eqref{phi}, \eqref{de1} and \eqref{dep}) that
\begin{align*}
\Phi(\lambda) = \frac{ \Phi^{\prime\prime}(0)  }{2} \lambda^2 +
\frac{ \Phi^{\prime\prime\prime}(0)  }{6} \lambda^3 + O_{\theta_0}( |z|^2 \lambda^4 ), \qquad \forall - \frac{\pi - \theta_0}{4} \leq \lambda \leq \frac{\pi - \theta_0}{4}, \  0 \leq \theta \leq \theta_0,
\end{align*}
with $\Phi^{\prime\prime\prime}(0) = O_{\theta_0}(|z|^2)$.
Furthermore, if $|z|^2 \delta^3$ is small enough (we will choose $\delta$ later to satisfy this condition) and $|z|$ large enough (cf. \eqref{ep11}), we can write
\begin{align*}
e^{ \Phi(\lambda) } &= e^{ \frac{ \Phi^{\prime\prime}(0)  }{2} \lambda^2 }
\left( 1 + \frac{ \Phi^{\prime\prime\prime}(0)  }{6} \lambda^3 + O_{\theta_0}( |z|^2 \lambda^4 ) + O_{\theta_0}( |z|^4 \lambda^6 )  \right),
\qquad \forall -\delta \leq \lambda \leq \delta,
\end{align*}
and combining it with \eqref{1h1}, we get
\begin{align*}
h(\lambda + \imath \theta)  e^{ \Phi(\lambda) } &=  e^{ \frac{ \Phi^{\prime\prime}(0)  }{2} \lambda^2 } h( \imath \theta ) \left( 1 + \frac{h^{\prime}( \imath \theta )}{h( \imath \theta )} \lambda +
\frac{\Phi^{\prime\prime\prime}(0)}{6} \lambda^3 + O_{\theta_0}( \lambda^2 + |z|^2 \lambda^4 + |z|^4 \lambda^6 )  \right)
\end{align*}
for all $-\delta \leq \lambda \leq \delta$.

Observe that
\begin{eqnarray*}
\int_{-\delta}^{\delta}  \lambda e^{ \frac{ \Phi^{\prime\prime}(0)  }{2} \lambda^2 } d\lambda =
\int_{-\delta}^{\delta}  \lambda^3 e^{ \frac{ \Phi^{\prime\prime}(0)  }{2} \lambda^2 } d\lambda =0
\end{eqnarray*}
and following from \eqref{dep}, we have
\begin{align*}
\int_{-\delta}^{\delta} |z|^{2j} \lambda^{2+2j} e^{ \frac{ \Phi^{\prime\prime}(0)  }{2} \lambda^2 } \, d\lambda \leq |z|^{2j} \int_{\R} \lambda^{2+2j}  e^{ \frac{ \Phi^{\prime\prime}(0)  }{2} \lambda^2 } d\lambda
&\sim |z|^{2j} \left( -\frac{\Phi^{\prime\prime}(0)}{2} \right)^{-j-\frac{3}{2}} \\
&\lesssim \left( -\frac{\Phi^{\prime\prime}(0)}{2} \right)^{-\frac{3}{2}},
\end{align*}
for $j= 0, 1, 2$. Hence, we yield from the above estimates and \eqref{ep11}
\begin{align*}
T_1 = h( \imath \theta)
\int_{-\delta}^{\delta} e^{ \frac{ \Phi^{\prime\prime}(0)}{2} \lambda^2 } \, d \lambda + O_{\theta_0}(d^{-3}(z, t)).
\end{align*}
Now note that
\begin{align*}
\int_{-\delta}^{\delta} e^{ \frac{ \Phi^{\prime\prime}(0)}{2} \lambda^2 } \, d \lambda &= \int_{\R} e^{ \frac{ \Phi^{\prime\prime}(0)  }{2} \lambda^2 } \, d \lambda
- \int_{|\lambda|> \delta} e^{ \frac{ \Phi^{\prime\prime}(0)  }{2} \lambda^2 } \, d \lambda \\
&= \sqrt{\pi} \left( -\frac{\Phi^{\prime\prime}(0)}{2} \right)^{-\frac{1}{2}} \left( 1 +
O\left( e^{  \frac{ \Phi^{\prime\prime}(0) }{4} \delta^2 }  \right) \right),
\end{align*}
and consequently,
\begin{align*}
T_1 = \sqrt{\pi} h( \imath \theta) \left( -\frac{\Phi^{\prime\prime}(0)}{2} \right)^{-\frac{1}{2}} \left( 1 +
O\left( e^{  \frac{ \Phi^{\prime\prime}(0) }{4} \delta^2 }  \right) \right) + O_{\theta_0}(d^{-3}(z, t)),
\end{align*}
under the assumption that $0 \leq \theta = \theta(z, t) \leq \theta_0$, $|z|^2 \delta^3 \ll 1$ and $|z|$ is large enough.

On the other hand, by Lemma \ref{l5} and \eqref{h1h2}, we have
\begin{eqnarray*}
\left| T_2 \right| \leq e^{-A_0 | z |^2 \delta^2} \int_{ |\lambda| > \delta } |h(\lambda + \imath \theta)| \, d\lambda  \lesssim_{\theta_0} e^{-A_0 | z |^2 \delta^2}, \quad \forall 0 \leq \theta \leq \theta_0.
\end{eqnarray*}

For $|z|$ large enough with $0 \leq \theta \leq \theta_0$, by setting
$$\delta = |z|^{-\frac{3}{4}} \ll 1,$$
we have
\begin{align*}
|z|^2 \delta^3 = |z|^{-\frac{1}{4}} \sim_{\theta_0} d^{-\frac{1}{4}}(z, t) \ll 1  \quad \mbox{and} \quad
- \Phi^{\prime\prime}(0) \delta^2 \sim_{\theta_0} |z|^2 \delta^2 = |z|^{\frac{1}{2}} \sim_{\theta_0} d^{\frac{1}{2}}(z, t).
\end{align*}
This completes the proof of Theorem \ref{t1}.
~ \hspace*{20pt} ~ \hfill $\Box$

\medskip

\renewcommand{\theequation}{\thesection.\arabic{equation}}
\section{Uniform asymptotic bounds of $M$ for $d(z, t) \longrightarrow +\infty$ with $0 < \epsilon = \pi - \theta \leq \epsilon_0 \ll 1$ and $D_1 + D_2 \gg 1$ when $ l\geq 2$ and $k_l = 1$}  \label{st2}
\setcounter{equation}{0}

Fix $0 < \epsilon_0 < \frac{1}{16} \min\left\{1, \frac{1 - a_{l -1}}{a_{l - 1}} \pi\right\}$, which is small enough such that (see \eqref{d2j})
\begin{align} \label{nd2j}
D_2 \sim |z|^2 \epsilon^2, \qquad \forall 0 \leq \epsilon \leq \epsilon_0.
\end{align}

Under our assumption,
setting in the sequel (see the definition of $h$ and $h_1$ in \eqref{defh} with $k_l = 1$),
\begin{align} \label{defs}
s(\xi) = s(\mathcal{K},\mathcal{A};\xi)  &= \frac{\xi}{\pi} h_1(\imath (\pi - \xi)) \prod_{j=1}^{l-1} h_1^{k_j}(\imath a_j (\pi - \xi)) \\
&= \frac{\xi}{\sin{\xi}}
\left( 1- \frac{\xi}{\pi} \right) \prod_{j=1}^{l-1} \left[ \frac{ a_j (\pi - \xi) }{ \sin(a_j (\pi - \xi) ) } \right]^{k_j},
\end{align}
then we have
\begin{align*}
h( \lambda + \imath \theta ) = \frac{\pi}{\epsilon + \imath \lambda} s( \epsilon + \imath \lambda ), \qquad \forall \lambda \in \R.
\end{align*}

Recall the definition of $M$ (see \eqref{spe} with $k_l = 1$) and by the change of variables $\lambda = \epsilon \tau$, we can write
\begin{align} \label{Md1}
M &= \epsilon \int_{\R} h(\imath (\pi - \epsilon (1 + \imath \tau))) e^{\Phi(\epsilon\tau)} \, d\tau = \pi \int_{\R} \frac{1}{1 + \imath \tau} s(\epsilon (1 + \imath \tau)) e^{\Phi(\epsilon \tau)} \, d\tau \nonumber\\
&= \pi \left( \int_{|\tau| \leq \delta_1} + \int_{ \delta_1 < |\tau| \leq \delta_2 }
+ \int_{|\tau| > \delta_2} \right) \nonumber \\
&= \pi (M_1 + M_2 + M_3),
\end{align}
where
\begin{align*}
\delta_2 = \frac{c_0}{\epsilon} \mbox{ with $c_0$ small enough, and } \delta_1 = \delta_1(\mathcal{K}, \mathcal{A}; D_1 + D_2) \ll 1,
\end{align*}
will be determined later.

It follows from Lemma \ref{l5} and \eqref{h1h2} that
\begin{align} \label{M3}
|M_3| \leq 2 \, e^{-\Re \Phi(c_0)} \int_{|\lambda| > c_0} |h(\lambda + \imath \theta)| \, d\lambda \leq
C(c_0, \epsilon_0) e^{-\Re \Phi(c_0)},
\end{align}
where the constant $C(c_0, \epsilon_0 ) > 0$ is independent of $(z, t)$.

\subsection{The estimate of $M_1$ (cf. \eqref{Md1})}

By the facts that the functions $\frac{1}{1 + \xi}$ and $s(\xi)$ (cf. \eqref{defs}) are holomorphic on the domain
\begin{align*}
\Omega = \left\{ \xi \in \C; |\xi| \leq \frac{3}{4} \min\left\{1, \frac{1 - a_{l -1}}{a_{l - 1}} \pi\right\} \right\},
\end{align*}
and $s(r) \geq \frac{1}{2}$ for all $0 \leq r \leq \frac{1}{16}$, we have
\begin{align}\label{lsl}
\frac{s(\epsilon (1 + \imath \tau))}{ 1 + \imath \tau } &= s(\epsilon) \left( 1 - \frac{s(\epsilon) - s'(\epsilon) \epsilon}{s(\epsilon)} \imath \tau + O(\tau^2) \right), \quad \forall 0 \leq \epsilon \leq \epsilon_0, \ -\epsilon_0 \leq \tau \leq \epsilon_0.
\end{align}

By the definition of $\Phi$ (cf. \eqref{p1p}), and $D_i$ ($i = 1, 2$) (see \eqref{ie1}), we get
\begin{align} \label{np1}
\Phi(\epsilon \tau ) = D_1 \left( \frac{1}{ 1 + \imath \tau } + \imath \tau - 1 \right) - D_2 \tau^2 + K(\epsilon \tau).
\end{align}

Note that $G$ (see \eqref{Ge2}) is holomorphic on $\Omega$, and that
\begin{align} \label{gk}
|G^{(k)}(\xi)| \leq C(\mathcal{K}, \mathcal{A}; k) |z|^2, \quad \forall |\xi| \leq 15 \epsilon_0, \  k = 2, 3, \cdots.
\end{align}
Using Taylor's formula for $G$ (so for $K$, cf. \eqref{ke}), we obtain
\begin{align} \label{pne}
\Phi(\epsilon \tau ) = -(D_1 + D_2)\tau^2 + \imath \left( D_1 - \frac{G^{\prime\prime\prime}(\epsilon)}{6} \epsilon^3 \right) \tau^3
+ O((D_1 + D_2)\tau^4), \qquad \forall - \epsilon_0 \leq \tau \leq \epsilon_0,
\end{align}
for all $(z, t)$ satisfying $t \geq 0$, $\epsilon = \epsilon(z, t) \leq \epsilon_0$ and $D_1 + D_2 \geq 1$.

In addition, when $D_1 + D_2 \geq \frac{1}{\epsilon_0}$ and $(D_1 + D_2) \delta_1^3 \leq 1$
(we will choose $\delta_1$ later to satisfy this condition), we yield from \eqref{pne} that
\begin{align*}
e^{\Phi(\epsilon \tau )} &= e^{-(D_1 + D_2)\tau^2}  \left( 1 + \imath \left( D_1 - \frac{G^{\prime\prime\prime}(\epsilon)}{6} \epsilon^3 \right) \tau^3
+ O((D_1 + D_2) \tau^4 +  (D_1 + D_2)^2 \tau^6) \right)
\end{align*}
for all $-\delta_1 \leq \tau \leq \delta_1$ and $0 \leq \epsilon \leq \epsilon_0$.
Combining this with \eqref{lsl}, we get
\begin{align*}
\frac{1}{ 1 + \imath \tau } s(\epsilon (1 + \imath \tau)) e^{\Phi(\epsilon \tau )} &= e^{-(D_1 + D_2)\tau^2} s(\epsilon)
\left[1 +  \frac{s^{\prime}(\epsilon) \epsilon - s(\epsilon)}{s(\epsilon)} \imath \tau
+ \imath  \left( D_1 - \frac{G^{\prime\prime\prime}(\epsilon)}{6} \epsilon^3 \right) \tau^3 \right] \\
 & \quad   + e^{-(D_1 + D_2)\tau^2} s(\epsilon) \left[  O(\tau^2 + (D_1 + D_2) \tau^4 + (D_1 + D_2)^2 \tau^6) \right],
\end{align*}
for all $-\delta_1 \leq \tau \leq \delta_1$ and $0 \leq \epsilon \leq \epsilon_0$.

Similar to the proof of Theorem \ref{t1}, we have
\begin{eqnarray*}
\int_{-\delta_1}^{\delta_1}  \tau e^{ -(D_1 + D_2) \tau^2 } d\tau =
\int_{-\delta_1}^{\delta_1}  \tau^3 e^{ -(D_1 + D_2) \tau^2 } d\tau  = 0,
\end{eqnarray*}
\begin{align*}
\left| \int_{-\delta_1}^{\delta_1} (D_1 + D_2)^{j} \tau^{2+2j} e^{ - (D_1 + D_2) \tau^2 } \, d\tau \right|
&\leq (D_1 + D_2)^{j} \int_{\R} \tau^{2+2j}  e^{ - (D_1 + D_2) \tau^2 } \, d\tau\\
&\sim (D_1 + D_2)^{j} (D_1 + D_2)^{- j - \frac{3}{2}}
\sim (D_1 + D_2)^{ - \frac{3}{2}},
\end{align*}
for $j= 0, 1, 2$ and
\begin{align*}
\int_{-\delta_1}^{\delta_1} e^{ - (D_1 + D_2) \tau^2 } \, d\tau = \int_{\R} - \int_{|\tau| \geq \delta_1} = \sqrt{\pi} (D_1 + D_2)^{ - \frac{1}{2}}
\left[ 1+ O\left(e^{-\frac{(D_1 + D_2)}{2} \delta_1^2}\right) \right].
\end{align*}

Consequently, under the assumption that $(D_1 + D_2) \delta_1^3 \leq 1$, we have
\begin{align} \label{M1e}
M_1 = s(\epsilon) \sqrt{\pi} (D_1 + D_2)^{ - \frac{1}{2}}
\left[ 1+ O\left(e^{-\frac{(D_1 + D_2)}{2} \delta_1^2}
+  (D_1 + D_2)^{ - 1} \right) \right],
\end{align}
for all $(z, t)$ satisfying $0 < \epsilon \leq \epsilon_0$ with $D_1 + D_2 \geq \frac{1}{\epsilon_0}$.

\subsection{The upper bound of $M_2$ (cf. \eqref{Md1})}  \label{S52}

Since $| \epsilon \tau | \leq c_0$ and $c_0$ is small, by \eqref{np1} and \eqref{gk} we have
\begin{align*}
\Phi(\epsilon \tau ) = D_1 \left( \frac{1}{ 1 + \imath \tau } + \imath \tau - 1 \right) - D_2 \tau^2 + O(|z|^2 \epsilon^3 \tau^3).
\end{align*}
Then for $c_0$ small enough, it deduces from the equation above and \eqref{nd2j} that
\begin{eqnarray}\label{RePhi1}
\Re\Phi(\epsilon \tau) \leq - \frac{D_2}{2} \tau^2, \qquad \forall \delta_1 \leq |\tau| \leq \delta_2 = \frac{c_0}{\epsilon}.
\end{eqnarray}

On the other hand, under the assumption that $(D_1 + D_2) \delta_1^3 \leq 1$ (as before, we will choose $\delta_1$ later) and $\delta_1 \leq \epsilon_0$, it follows from Lemma \ref{l5} and \eqref{pne} that
\begin{eqnarray}\label{RePhi2}
\Re\Phi(\epsilon \tau) \leq \Re\Phi(\epsilon \delta_1) = -(D_1 + D_2)\delta_1^2 + O(1), \quad \forall |\tau| > \delta_1.
\end{eqnarray}
Taking average of \eqref{RePhi1} and \eqref{RePhi2}, we have
\begin{eqnarray} \label{RePhi3}
\Re\Phi(\epsilon \tau) \leq  - \frac{D_2}{4} \tau^2 - \frac{(D_1 + D_2)}{2} \delta_1^2 + O(1), \quad
\forall \delta_1 \leq |\tau| \leq \delta_2 = \frac{c_0}{\epsilon}.
\end{eqnarray}
We divide it into cases.

Case (1): $D_2 \geq 1$. If $(D_1 + D_2) \delta_1^3 \leq 1$ and $\delta_1 \leq \epsilon_0$, then from \eqref{RePhi3} we have
\begin{align} \label{M2a}
|M_2| \lesssim e^{ - \frac{(D_1 + D_2)}{2} \delta_1^2} \int_{\delta_1 < |\tau| \leq \delta_2 } e^{ - \frac{D_2 }{4} \tau^2} \, d\tau \lesssim \frac{1}{\sqrt{D_2}} e^{ - \frac{(D_1 + D_2)}{2} \delta_1^2} = O\left( e^{ - \frac{(D_1 + D_2)}{2} \delta_1^2} \right).
\end{align}

Case (2): $D_2 < 1$. We have $D_1 \geq 1$ when $D_1 + D_2 \geq \frac{1}{\epsilon_0} \geq 2$. We set
\begin{align}\label{defH}
\qquad H(u) = H(\mathcal{K},\mathcal{A};\epsilon;u) = s (\epsilon + \imath u),
\end{align}
and write
\begin{eqnarray*}
M_2 = \int_{\delta_1}^{\delta_2} + \int_{-\delta_2}^{-\delta_1}  = M_{2,1} + M_{2,2}.
\end{eqnarray*}
Note that we only need to estimate $M_{2,1}$ since the estimate of the other term is similar. Recalling the definition of $M_2$, cf. \eqref{Md1}, and by \eqref{np1}, we apply integration by parts to obtain that
\begin{align*}
M_{2,1} &= \int_{\delta_1}^{\delta_2} \frac{H(\epsilon \tau)}{1 + \imath \tau}
e^{- D_1 + \imath D_1 \tau - D_2 \tau^2 + \frac{D_1}{1 + \imath \tau} + K(\epsilon \tau) } d\tau \\
&= \int_{\delta_1}^{\delta_2} \frac{H(\epsilon \tau)}{( 1 + \imath \tau ) ( -2 D_2 \tau + \imath D_1 )}
e^{ - D_1 + \frac{D_1}{1 + \imath \tau} + K(\epsilon \tau) }
\frac{d}{d\tau} \left( e^{ \imath D_1 \tau - D_2 \tau^2 } \right) d\tau\\
&= \sum_{m = 1}^4 M_{2, 1, m},
\end{align*}
with
\begin{align*}
M_{2, 1, 1} =  \frac{H(\epsilon \delta_2)}{( 1 + \imath \delta_2 ) ( -2 D_2 \delta_2 + \imath D_1 )} e^{\Phi(\epsilon \delta_2)}, \quad
M_{2, 1, 2} =  - \frac{H(\epsilon \delta_1)}{( 1 + \imath \delta_1 ) ( -2 D_2 \delta_1 + \imath D_1 )} e^{\Phi(\epsilon \delta_1)},
\end{align*}
and
\begin{align*}
M_{2, 1, 3} &=  - \int_{\delta_1}^{\delta_2} e^{\Phi(\epsilon \tau)}
\frac{d}{d\tau} \Big[ \frac{H(\epsilon \tau)}{( 1 + \imath \tau ) ( -2 D_2 \tau + \imath D_1 )} \Big] \, d\tau, \\
M_{2, 1, 4} &= - \int_{\delta_1}^{\delta_2} \frac{H(\epsilon \tau) e^{\Phi(\epsilon \tau)}}{( 1 + \imath \tau ) ( -2 D_2 \tau + \imath D_1 )} \frac{d}{d\tau} \Big[ - D_1 + \frac{D_1}{1 + \imath \tau} + K(\epsilon \tau) \Big] \, d\tau.
\end{align*}

Now, by $D_1 \geq 1$ and the fact that $H(\tau) = O(1)$
for all $|\tau| \leq c_0$ with $0 \leq \epsilon \leq \epsilon_0$ (cf. \eqref{defH} and \eqref{defs}), \eqref{RePhi2} implies that
\begin{align*}
|M_{2, 1, 1}| + |M_{2, 1, 2}| = O\left(e^{-(D_1 + D_2) \delta_1^2}\right).
\end{align*}

Note that
\begin{align*}
M_{2, 1, 3} &= - \int_{\delta_1}^{\delta_2} e^{\Phi(\epsilon \tau)}
\frac{\epsilon H^{\prime}(\epsilon \tau)}{( 1 + \imath \tau ) ( -2 D_2 \tau + \imath D_1 )} \, d\tau
+ \imath \int_{\delta_1}^{\delta_2} e^{\Phi(\epsilon \tau)}
\frac{H(\epsilon \tau)}{( 1 + \imath \tau )^2 ( -2 D_2 \tau + \imath D_1 )} \, d\tau \\
&+ \int_{\delta_1}^{\delta_2} e^{\Phi(\epsilon \tau)}
\frac{( -2 D_2 ) H(\epsilon \tau)}{( 1 + \imath \tau ) ( -2 D_2 \tau + \imath D_1 )^2} \, d\tau = M_{2, 1, 3, 1} + M_{2, 1, 3, 2} + M_{2, 1, 3, 3}.
\end{align*}
Similarly, from $D_1 \geq 1$, $H(\tau) = O(1)$ for all $|\tau| \leq c_0$ with $0 \leq \epsilon \leq \epsilon_0$ and \eqref{RePhi2} we have
\begin{align*}
|M_{2, 1, 3, 2}| &\leq e^{\Re \Phi(\epsilon \delta_1)} \sup_{0 \leq s \leq c_0} |H(s)| \int_{\delta_1}^{\delta_2} \frac{1}{1 + \tau^2} \, d\tau = O\left(e^{-(D_1 + D_2) \delta_1^2}\right), \\
|M_{2, 1, 3, 3}| &\leq e^{\Re \Phi(\epsilon \delta_1)} \sup_{0 \leq s \leq c_0} |H(s)| \int_{\delta_1}^{\delta_2} \frac{2 D_2}{\sqrt{(1 + \tau^2)(1 + (D_2 \tau)^2)}} \, d\tau = O\left(e^{-(D_1 + D_2) \delta_1^2}\right),
\end{align*}
and by the additional fact that $H'(\tau) = O(1)$ for all $0 \leq \epsilon \leq \epsilon_0$ and $|\tau| \leq c_0$, we have
\begin{align*}
|M_{2, 1, 3, 1}| &\leq e^{\Re \Phi(\epsilon \delta_1)} \sup_{0 \leq s \leq c_0} |H'(s)| \int_{\delta_1}^{\delta_2} \epsilon \, d\tau = O\left(e^{-(D_1 + D_2) \delta_1^2}\right).
\end{align*}

Moreover, observe that
\begin{align*}
M_{2, 1, 4} &= \imath \int_{\delta_1}^{\delta_2} e^{\Phi(\epsilon \tau)}
\frac{D_1 H(\epsilon \tau)}{( 1 + \imath \tau )^3 ( -2 D_2 \tau + \imath D_1 )} \, d\tau  - \int_{\delta_1}^{\delta_2} e^{\Phi(\epsilon \tau)} \frac{\epsilon H(\epsilon \tau) K^{\prime}(\epsilon \tau)}{( 1 + \imath \tau ) ( -2 D_2 \tau + \imath D_1 )} \, d\tau \\
&= M_{2, 1, 4, 1} + M_{2, 1, 4, 2},
\end{align*}
and $M_{2, 1, 4, 1} = O\left(e^{-(D_1 + D_2) \delta_1^2}\right)$ from a similar argument used above. It remains to establish the estimate of $M_{2, 1, 4, 2}$. It deduces from \eqref{ke} and \eqref{gk} that
\begin{align*}
K^{\prime}(u) = O(|z|^2 u^2), \qquad \forall 0 \leq \epsilon \leq \epsilon_0, \ |u| \leq c_0,
\end{align*}
and combining it with \eqref{RePhi3} and \eqref{nd2j}, we obtain from the change of variables $u = \epsilon \tau$ that
\begin{eqnarray*}
|M_{2,1, 4,2}| \lesssim  e^{ - \frac{(D_1 + D_2)}{2} \delta_1^2} \int_{\epsilon \delta_1 }^{c_0} e^{-c |z|^2 u^2} |z|^2 u^2 \, du \lesssim e^{ - \frac{(D_1 + D_2)}{2} \delta_1^2},
\end{eqnarray*}
where $c$ is the implicit constant in the estimate \eqref{nd2j}.

Consequently, under the assumption that $(D_1 + D_2) \delta_1^3 \leq 1$, $\delta_1 \leq \epsilon_0$, $D_1 + D_2 \geq \frac{1}{\epsilon_0}$ and $D_2 < 1$, we have
\begin{align} \label{M2t}
M_2 = O\left(e^{ - \frac{(D_1 + D_2)}{2} \delta_1^2}\right).
\end{align}

When $0 < \epsilon \leq \epsilon_0$ and $D_1 + D_2 \ge \epsilon_0^{-3} > \epsilon_0^{-1}$,
take
$$\delta_1 =  (D_1 + D_2)^{-\frac{1}{3}} \leq \epsilon_0 \ll 1,$$
and we get
\begin{align*}
(D_1 + D_2) \delta_1^4 = \delta_1 \ll 1, \quad (D_1 + D_2) \delta_1^3 = 1  \quad \mbox{and} \quad (D_1 + D_2) \delta_1^2 = (D_1 + D_2)^{\frac{1}{3}} \gg 1.
\end{align*}
In conclusion, by \eqref{Md1}, \eqref{M3}, \eqref{M1e}, \eqref{M2a}, \eqref{M2t} and the fact that
$$\Re \Phi(c_0) \leq \Re \Phi(\epsilon \delta_1) \leq -\frac{D_1 + D_2}{2} \delta_1^2,$$
we have
\begin{align*}
M &= \pi^{\frac{3}{2}} s(\epsilon) \, (D_1 + D_2)^{ - \frac{1}{2}} \left( 1 + O\left( (D_1 + D_2)^{ - 1} \right) \right) \\
&= \pi^{\frac{3}{2}} s(0) \, (D_1 + D_2)^{ - \frac{1}{2}}(1 + O(\epsilon))\left( 1 + O\left( (D_1 + D_2)^{ - 1} \right) \right).
\end{align*}

This completes the proof of Theorem \ref{t2}.
\hspace*{20pt} ~ \hfill $\Box$

\medskip

\renewcommand{\theequation}{\thesection.\arabic{equation}}
\section{Uniform asymptotic estimates of $M$ for $d(z, t) \longrightarrow +\infty$ with $0 \leq \epsilon = \pi - \theta \leq \epsilon_0$ and $D_1 + D_2 \leq \gamma_0$ when $l \geq 2$ and $k_l = 1$}  \label{st3}
\setcounter{equation}{0}

Let $\epsilon_0$ be defined as in the beginning of Section \ref{st2}.

Under our assumption,
we split the second integral in \eqref{Md1} into two pieces as follows
\begin{align} \label{ndm}
M = \pi \int_{\R} \frac{1}{1 + \imath \tau} s(\epsilon (1 + \imath \tau)) e^{\Phi(\epsilon \tau)} \, d\tau
= \pi \left( \int_{|\tau| \leq \tau_0} + \int_{  |\tau| > \tau_0 } \right)= \pi (M_1^* + M_2^* ),
\end{align}
where $\tau_0 = \tau_0(\mathcal{K},\mathcal{A};z , t) \gg 1$ will be determined later.

When $\epsilon \tau_0$ small enough, by the fact that $s$ (cf. \eqref{defs}) and $K$ (cf. \eqref{ke}) are holomorphic on $\{\xi \in \C; |\xi| < 8 \epsilon_0 \}$ and $s(0) \geq 1$, using the mean value theorem we have
\begin{align*}
s(\epsilon (1 + \imath \tau)) = s(0) \left( 1  + O(\epsilon |1 + \imath \tau|) \right) = s(0) \left( 1  + O(\epsilon \sqrt{1 + \tau^2}) \right),
\end{align*}
and
\begin{align*}
K(\epsilon \tau) = K(\imath \epsilon) + \left( K(\epsilon \tau) - K(\imath \epsilon) \right) = J_* + O(|z|^2 \epsilon^3 |1 + \imath \tau|) =  J_*  + O(D_2 \epsilon \sqrt{1 + \tau^2}),
\end{align*}
for all $0 \leq \epsilon \leq \epsilon_0$ and $|\epsilon\tau| \le \epsilon \tau_0 \leq 5 \epsilon_0$,
where we have used \eqref{nd2j}
in the last equality. Combining these two expansions we obtain
\begin{align*}
s(\epsilon (1 + \imath \tau)) e^{K(\epsilon \tau)} = s(0) \, e^{J_*} \left(1 +  O_{\gamma_0}(\epsilon \sqrt{1 + \tau^2})
\right), \qquad \forall 0 \leq \epsilon \leq \epsilon_0, \ |\epsilon\tau| \le \epsilon \tau_0 \leq 5 \epsilon_0.
\end{align*}

Set
\begin{eqnarray*}
S(\tau) = \imath D_1 \tau - D_2 \tau^2 + \frac{D_1}{1 + \imath \tau}.
\end{eqnarray*}

\subsection{The estimate of $M_1^*$ (cf. \eqref{ndm})}

By \eqref{np1} and the expansion above, we have $\Phi(\epsilon \tau) = - D_1 + S(\tau) +K(\epsilon \tau)$ and
\begin{align} \label{M*1}
M_1^* &= s(0) e^{- D_1 + J_*} \left[ \int_{|\tau| \leq \tau_0} \frac{1}{1 + \imath \tau} e^{S(\tau)} \, d\tau
+ \int_{|\tau| \leq \tau_0} O_{\gamma_0}\left(\epsilon e^{\Re S(\tau)}\right) \, d\tau \right] \nonumber \\
&= s(0) e^{- D_1 + J_*} \left[ \int_{\R} \frac{1}{1 + \imath \tau} e^{S(\tau)} \, d\tau
-\int_{|\tau| > \tau_0} \frac{1}{1 + \imath \tau} e^{S(\tau)} \, d\tau
+ \int_{|\tau| \leq \tau_0} O_{\gamma_0}\left(\epsilon e^{\Re S(\tau)}\right) \, d\tau \right] \nonumber \\
&= s(0) e^{- D_1 + J_*} (\bar{S}_1 - S_1^* + S_0^*).
\end{align}

Recalling the definition of $S_1$ (cf. \eqref{SKL}) and using Proposition \ref{prop1} in Subsection \ref{precise} (with $\nu = 1$, $r = D_1$ and $b = D_2$), we get
\begin{align} \label{M*1a}
\bar{S}_1 =   S_1 \sim_{\gamma_0} 1
\end{align}
and a direct computation yields that
\begin{align} \label{M*1b}
|S_0^*| =  e^{D_1} \epsilon \tau_0 \, O_{\gamma_0}(1) =\epsilon \tau_0 \, O_{\gamma_0}(1).
\end{align}

Now we consider the term $S_1^*$. In fact, the estimate can be found in \cite[Cas 2, p. 816]{Li12} and we include it for the sake of completeness. We can write
\begin{eqnarray*}
S_1^* = \int_{\tau_0}^{+\infty} + \int_{-\infty}^{-\tau_0} = S_{1,1}^* + S_{1,2}^*.
\end{eqnarray*}
As before, we only estimate $S_{1,1}^*$ and the estimate for the other term is similar.

Observe that
\begin{eqnarray*}
S_{1,1}^* = \int_{\tau_0}^{+\infty} \frac{e^{\frac{D_1}{1 + \imath \tau}}}{( 1 + \imath \tau ) ( -2 D_2 \tau + \imath D_1 )} \,
\frac{d}{d\tau}\left( e^{ \imath D_1 \tau - D_2 \tau^2 } \right) d\tau = \sum_{m=1}^4 S_{1,1,m}^*,
\end{eqnarray*}
with
\begin{align*}
S_{1,1,1}^* &= - \frac{e^{\frac{D_1}{1 + \imath \tau_0}}}{( 1 + \imath \tau_0 ) ( -2 D_2 \tau_0
+ \imath D_1 )} e^{ \imath D_1 \tau_0 - D_2 \tau_0^2 }, \\
S_{1,1,2}^* &= \imath \int_{\tau_0}^{+\infty} \frac{e^{\frac{D_1}{1 + \imath \tau}}}{( 1 + \imath \tau )^2 ( -2 D_2 \tau + \imath D_1 )} e^{ \imath D_1 \tau - D_2 \tau^2 } \, d\tau, \\
S_{1,1,3}^* &= \imath D_1 \int_{\tau_0}^{+\infty} \frac{ e^{\frac{D_1}{1 + \imath \tau}}}{( 1 + \imath \tau )^3 ( -2 D_2 \tau + \imath D_1 )} e^{ \imath D_1 \tau - D_2 \tau^2 } \, d\tau,
\end{align*}
and
\begin{align*}
S_{1,1,4}^* = - 2 D_2 \int_{\tau_0}^{+\infty} \frac{ e^{\frac{D_1}{1 + \imath \tau}}}{( 1 + \imath \tau ) ( -2 D_2 \tau + \imath D_1 )^2}
e^{ \imath D_1 \tau - D_2 \tau^2 } \, d\tau.
\end{align*}
It is not hard to find that
\begin{align*}
\sum_{m = 1}^4 |S_{1,1, m}^*| =  \frac{e^{\frac{D_1}{1 + \tau_0^2} -  D_2 \tau_0^2}}{D_1 \tau_0 + D_2 \tau_0^2} \, O\left( 1 +     \frac{D_1}{\tau_0} \right).
\end{align*}

Consequently,
\begin{align} \label{S1*}
S_1^* =  \frac{e^{\frac{D_1}{1 + \tau_0^2} - D_2 \tau_0^2}}{D_1 \tau_0 + D_2 \tau_0^2} \, O\left( 1 +     \frac{D_1}{\tau_0}  \right).
\end{align}

\subsection{The upper bound of $M_2^*$ (cf. \eqref{ndm})} \label{S52*}

We write
\begin{eqnarray*}
M_2^* = \int_{\tau_0}^{+\infty} + \int_{-\infty}^{-\tau_0} = M_{2,1}^* + M_{2,2}^*,
\end{eqnarray*}
and it suffices to estimate $M_{2,1}^*$. We argue as in Case (2) of Subsection \ref{S52}, and get
\begin{align*}
M_{2, 1}^* = M_{2, 1, 1}^* + M_{2, 1, 3, 1}^* + M_{2, 1, 3, 2}^* + M_{2, 1, 3, 3}^* + M_{2, 1, 4, 1}^* + M_{2, 1, 4, 2}^*,
\end{align*}
where
\begin{align*}
M_{2, 1, 1}^* &= - \frac{H(\epsilon \tau_0)}{( 1 + \imath \tau_0 ) ( -2 D_2 \tau_0 + \imath D_1 )} e^{\Phi(\epsilon \tau_0)},\\
M_{2, 1, 3, 1}^* &=  - \int_{\tau_0}^{+\infty} e^{\Phi(\epsilon \tau)}
\frac{\epsilon H^{\prime}(\epsilon \tau)}{( 1 + \imath \tau ) ( -2 D_2 \tau + \imath D_1 )} \, d\tau, \\
M_{2, 1, 3, 2}^* &= \imath \int_{\tau_0}^{+\infty} e^{\Phi(\epsilon \tau)}
\frac{H(\epsilon \tau)}{( 1 + \imath \tau )^2 ( -2 D_2 \tau + \imath D_1 )} \, d\tau, \\
M_{2, 1, 3, 3}^* &= \int_{\tau_0}^{+\infty} e^{\Phi(\epsilon \tau)}
\frac{( -2 D_2 ) H(\epsilon \tau)}{( 1 + \imath \tau ) ( -2 D_2 \tau + \imath D_1 )^2} \, d\tau,
\end{align*}
and
\begin{align*}
M_{2, 1, 4, 1}^* &= \imath \int_{\tau_0}^{+\infty} e^{\Phi(\epsilon \tau)}
\frac{D_1 H(\epsilon \tau)}{( 1 + \imath \tau )^3 ( -2 D_2 \tau + \imath D_1 )} \, d\tau, \\
M_{2, 1, 4, 2}^* &= - \int_{\tau_0}^{+\infty} e^{\Phi(\epsilon \tau)}
\frac{\epsilon H(\epsilon \tau) K^{\prime}(\epsilon \tau)}{( 1 + \imath \tau ) ( -2 D_2 \tau + \imath D_1 )} \, d\tau.
\end{align*}

Now Lemma \ref{l5} implies that $|e^{\Phi(\epsilon \tau)}| \leq 1$, and it follows from the definition of $H(u) =  s(\epsilon + \imath u)$
(cf. \eqref{defH} and \eqref{defs}) and \eqref{h1h2} that
\begin{align}
\label{estH}\int_{\R} \max\{u^2,1\} \left[ |H'(u)| + |H(u)| \right]  \, du = O(1), & \\
\label{estH*} |H(u)| = O(1), \qquad \forall 0 \leq \epsilon \leq \epsilon_0, u \in \R. &
\end{align}

As a result, it is not hard to see
\begin{align*}
|M^*_{2, 1, 1}| + |M^*_{2, 1, 3, 1}| + |M^*_{2, 1, 3, 2}| + |M^*_{2, 1, 4, 1}| = \frac{1}{D_1 \tau_0 + D_2 \tau_0^2} \, O\left(  1 + \frac{D_1}{\tau_0} \right),
\end{align*}
and
\begin{align*}
|M^*_{2, 1, 3, 3}| \lesssim \int_{\tau_0}^{+\infty} \frac{D_2}{\tau (D_1 + D_2 \tau)^2} \, d\tau = \frac{1}{D_1 \tau_0 + D_2 \tau_0^2}.
\end{align*}

It remains to estimate $|M^*_{2, 1, 4, 2}|$. Note that Lemma \ref{l5} implies that
\begin{align*}
\Re\Phi(u) \leq -A_0 |z|^2 \min\{u^2,1\}, \qquad \forall u \geq 0,
\end{align*}
and it follows from \eqref{ke} and \eqref{gk} (see also \eqref{Ge2}) that
\begin{align*}
|K'(u)| = |z|^2u^2O(1), \qquad \forall u \geq 0, \ 0 \leq \epsilon \leq \epsilon_0.
\end{align*}
Combining the two estimates above, we have
\begin{align*}
|K'(u)|e^{\Re\Phi(u)} \lesssim \max\{u^2,1\},
\end{align*}
which yields
\begin{eqnarray*}
|M_{2, 1, 4, 2}^*| \lesssim \frac{1}{D_1 \tau_0 + D_2 \tau_0^2} \int_0^{+\infty} \max\{u^2,1\} |H(u)| \, du
\lesssim \frac{1}{D_1 \tau_0 + D_2 \tau_0^2},
\end{eqnarray*}
where we have used $|e^{\Phi(\epsilon \tau)}| \leq 1$ in the first ``$\lesssim$'' and \eqref{estH} in the second ``$\lesssim$''.

Consequently, we have
\begin{align} \label{M2*e}
|M_2^*| = \frac{1}{D_1 \tau_0 + D_2 \tau_0^2} \, O\left(  1 + \frac{D_1}{\tau_0} \right).
\end{align}

For $(z,t)$ satisfying $d(z, t) \gg 1$ and $0 < \epsilon \leq \epsilon_0$, we observe that the first equality of \eqref{dd2} implies that
\begin{align} \label{d2eq}
d^2(z, t)  \sim \frac{|z_l|^2}{\epsilon^2} + |z|^2.
\end{align}

Setting
$$\tau_0 = d^{-\frac{3}{4}}(z, t) \epsilon^{-1},$$
we obtain
$$\epsilon \tau_0 = d^{-\frac{3}{4}}(z, t) \ll 1$$
and
\begin{align*}
D_1 \tau_0 + D_2 \tau_0^2 \sim \frac{|z_l|^2}{\epsilon^2} d^{-\frac{3}{4}}(z, t) + |z|^2 d^{-\frac{3}{2}}(z, t) \geq \left( \frac{|z_l|^2}{\epsilon^2} + |z|^2 \right) d^{-\frac{3}{2}}(z, t) \sim d^{\frac{1}{2}}(z, t) \gg 1,
\end{align*}
where the first ``$\sim$'' follows from \eqref{nd2j} and \eqref{ie1}, and the second ``$\sim$'' follows from \eqref{d2eq}.
Moreover, if $D_1 + D_2 \leq \gamma_0 (\ge 1)$, then we proceed in a similar way to obtain
\begin{align*}
2 \gamma_0 \tau_0 &\geq D_1 \tau_0 + \sqrt{D_2} \tau_0 \sim d^{-\frac{3}{4}}(z, t) \max\left\{\frac{|z_l|^2}{\epsilon^2}, |z|\right\} \\
&\geq d^{-\frac{3}{4}}(z, t) \max\left\{\frac{|z_l|}{\epsilon}, |z|\right\} \sim d^{\frac{1}{4}}(z, t) \gg 1,
\end{align*}
which justifies $\tau_0$ is big enough for fixed $\gamma_0$.

In conclusion, we get from \eqref{ndm}--\eqref{S1*} and \eqref{M2*e} that
\begin{align*}
M =\pi s(0) e^{- D_1 + J_*} S_1 \, \left(1 + O_{\gamma_0}\left(d^{-\frac{1}{2}}(z, t)\right)\right),
\end{align*}
under the assumption that $d(z, t) \gg 1$, $0 < \epsilon \leq \epsilon_0$ and $D_1 + D_2 \leq \gamma_0$.

This completes the proof of Theorem \ref{t3}.
~ \hspace*{20pt} ~ \hfill $\Box$

\medskip

\renewcommand{\theequation}{\thesection.\arabic{equation}}
\section{Proof of Theorem \ref{t4}}  \label{st4}
\setcounter{equation}{0}

With the simplified notation (replacing $p^{\H(\mathcal{K}, \mathcal{A})}$ and $\nabla_{\H(\mathcal{K}, \mathcal{A})}$ with $p$ and $\nabla$ respectively) in Theorem \ref{t4}, by induction on $|\beta|$,
we observe that $ \U_1^{\beta_1} \cdots \U_m^{\beta_m} p(g)$ is a linear combination of terms
of the following type (recall the definition of $\varphi$, $h$ and $h_1$ in \eqref{defvarphi} and \eqref{defh} respectively)
\begin{eqnarray*}
I_{\alpha,\iota} = P_{\alpha,\iota}(z)
\int_{\R}  h(\lambda) e^{\varphi(\lambda)}
\prod_{j=1}^l \left[ h_1^{\iota_j}(a_j \lambda)\sinh^{\iota_{j,1}}(a_j\lambda)\cosh^{\iota_{j,2}}(a_j\lambda) \right] \, d\lambda
\end{eqnarray*}
where $\alpha$, $\iota \in \N^l$ satisfies
\begin{align}\label{index}
|\alpha| = |\beta|, \qquad  \frac{\alpha_i}{2} \leq \iota_i = \iota_{i,1}+\iota_{i,2} \leq \alpha_i, \quad \forall 1 \leq i \leq l,
\end{align}
and $P_{\alpha,\iota}(z)$ is a homogeneous polynomial of degree $2 |\iota | - |\alpha|$  with complex coefficients as follows:
\begin{align*}
P_{\alpha,\iota}(z) = \sum_{\sum\limits_{j = 1}^{k_i} (q_{i, j} + r_{i, j}) = 2 \iota_i - \alpha_i, 1 \leq i \leq l} C(q_{i, j}, r_{i, j}) \prod_{i, j} x_{i, j}^{q_{i, j}} y_{i, j}^{r_{i, j}}.
\end{align*}
One typical example is
\begin{align*}
\X_{1,1}\Y_{1,1}p(g) &= -\frac{\imath}{2} \int_{\R} h(\lambda)e^{\varphi(\lambda)} h_1(a_1 \lambda)  \sinh(a_1\lambda) \, d\lambda \\
& + \frac{1}{4} x_{1,1} y_{1,1}\int_{\R} h(\lambda)e^{\varphi(\lambda)}
h_1^2(a_1 \lambda)\cosh^2(a_1\lambda) \, d\lambda \\
& + \frac{\imath}{4} (x_{1,1}^2 - y_{1,1}^2 ) \int_{\R} h(\lambda)e^{\varphi(\lambda)}
h_1^2(a_1 \lambda) \sinh(a_1\lambda) \cosh(a_1\lambda) \, d\lambda \\
& +  \frac{1}{4} x_{1,1} y_{1,1} \int_{\R} h(\lambda)e^{\varphi(\lambda)}
h_1^2(a_1 \lambda)\sinh^2(a_1\lambda) \, d\lambda.
\end{align*}

By the smoothness of $p$, without loss of generality, we may assume that $d(z, t)$ is large enough.  For fixed $\alpha$ and $\iota$, it follows from \eqref{l1} that
\begin{align*}
|P_{\alpha,\iota}(z) | \lesssim d(z, t)^{\sum\limits_{i = 1}^{l - 1} (2 \iota_i - \alpha_i)} |z_l|^{2 \iota_l - \alpha_l} \leq d(z, t)^{\sum\limits_{i = 1}^{l - 1} \alpha_i } |z_l|^{2 \iota_l - \alpha_l},
\end{align*}
where in the last inequality we have used $\iota_i \leq \alpha_i$ in \eqref{index}.
Combining it with Proposition \ref{prop1} in Subsection \ref{precise} (with $r = D_1$ and $b = D_2$ as before), we can show that
\begin{align} \nonumber
|I_{\alpha,\iota}| &\lesssim_{\epsilon_0} d(z, t)^{\sum\limits_{i = 1}^{l - 1}  \alpha_i } e^{-\frac{d^2(z,t)}{4}} \\
&\times |z_l|^{2 \iota_l - \alpha_l} \left\{\begin{array}{ll}
|z|^{-1} & \textrm{if $|\theta| \leq \pi -\epsilon_0$,} \\
\label{IAI}
(D_1+D_2)^{-\frac{1}{2}} \epsilon^{1 - k_l - \iota_l} & \textrm{if $\epsilon \leq \epsilon_0$ and $D_1 + D_2 \geq \epsilon_0^{-3}$,}\\
\left(\frac{D_1+\sqrt{D_2}}{\epsilon}\right)^{k_l+\iota_l-1} & \textrm{if $\epsilon \leq \epsilon_0$ and $D_1+D_2 \leq \epsilon_0^{-3}$.}\\
\end{array} \right.
\end{align}

In the case that $d(z, t) \gg 1$ and $|\theta| \leq \pi - \epsilon_0$, by \eqref{pehk}, \eqref{IAI} and \eqref{l1}, we get
$$\frac{|I_{\alpha,\iota}|}{p(z,t)} \lesssim_{\epsilon_0} d(z, t)^{(2 \iota_l - \alpha_l) + \sum\limits_{i = 1}^{l - 1}
\alpha_i } \leq d(z, t)^{\alpha_l + \sum\limits_{i = 1}^{l - 1} \alpha_i } = d(z, t)^{|\beta|},$$
where we have used the fact that $  \iota_l \leq \alpha_l$ and $|\alpha| = |\beta|$ (cf. \eqref{index}).

Now we consider the case where $d(z, t) \gg 1$ and $0 < \epsilon \leq \epsilon_0$.
Recall that (see \eqref{d2eq}, \eqref{nd2j} and \eqref{ie1})
\begin{align*}
d^2(z, t)  \sim \frac{|z_l|^2}{\epsilon^2} + |z|^2, \quad D_2 \sim |z|^2 \epsilon^2, \quad D_1 = \frac{\pi}{4} \frac{|z_l|^2}{\epsilon},
\end{align*}
which yields
\begin{align}\label{lbofd}
\epsilon d^2(z, t) \gtrsim D_1 + \frac{D_2}{\epsilon}.
\end{align}
We divide it into cases.

 Case (1): $D_1+D_2  \geq \epsilon_0^{-3}$.
From \eqref{pehk} and \eqref{IAI} we have
\begin{align*}
\frac{|I_{\alpha,\iota}|}{d(z, t)^{|\beta|} p(z,t)} \lesssim_{\epsilon_0} \frac{|z_l|^{2 \iota_l - \alpha_l}}{d(z, t)^{\alpha_l} \epsilon^{\iota_l}} &= \left( \frac{|z_l|}{\epsilon} \right)^{2 \iota_l - \alpha_l} \left(\epsilon d^2(z, t) \right)^{-(\alpha_l - \iota_l)} d(z, t)^{-(2 \iota_l - \alpha_l)} \\
&\lesssim d(z, t)^{2 \iota_l - \alpha_l} \left(\epsilon d^2(z, t) \right)^{-(\alpha_l - \iota_l)} d(z, t)^{-(2 \iota_l - \alpha_l)} \\
&\lesssim 1,
\end{align*}
where the first ``$\lesssim$'' follows from  \eqref{d2eq} and the second ``$\lesssim$'' follows from $  \iota_l \leq \alpha_l$ in \eqref{index} and $\epsilon d^2(z, t)  \gtrsim 1$ by \eqref{lbofd}.

Case (2): $D_1+D_2\leq \epsilon_0^{-3}$.
As before, \eqref{pehk} and \eqref{IAI} yield
\begin{align*}
\frac{|I_{\alpha,\iota}|}{d(z, t)^{|\beta|} p(z,t)} &\lesssim_{\epsilon_0} \left( \frac{D_1 + \sqrt{D_2}}{\epsilon} \right)^{\iota_l} \left( \frac{|z_l|}{\epsilon} \right)^{2 \iota_l - \alpha_l} \epsilon^{2 \iota_l - \alpha_l} d(z, t)^{-\alpha_l} \\
&\lesssim \left( \frac{D_1 + \sqrt{D_2}}{\epsilon} \right)^{\iota_l} d(z, t)^{2 \iota_l - \alpha_l} \epsilon^{2 \iota_l - \alpha_l} d(z, t)^{-\alpha_l} \\
&= \left( \frac{\epsilon d(z, t)^2}{D_1 + \sqrt{D_2}} \right)^{-(\alpha_l - \iota_l)} (D_1 + \sqrt{D_2})^{2 \iota_l - \alpha_l} \\
&\lesssim_{\epsilon_0} \left( \frac{\epsilon d(z, t)^2}{D_1 + \sqrt{D_2}} \right)^{-(\alpha_l - \iota_l)},
\end{align*}
where we have used \eqref{d2eq} for the ``$\lesssim$'' and $  \iota_l \leq \alpha_l$ in  \eqref{index} for the second ``$\lesssim_{\epsilon_0}$''.
It remains to show that $\frac{\epsilon d(z, t)^2}{D_1 + \sqrt{D_2}} \gtrsim 1$. Indeed, using \eqref{lbofd}, \eqref{nd2j} and \eqref{d2eq}, we can show
\begin{eqnarray*}
\frac{\epsilon d^2(z, t)}{D_1 + \sqrt{D_2}} \gtrsim \frac{D_1 + \frac{D_2}{\epsilon}}{D_1 + \sqrt{D_2}}\gtrsim \left\{\begin{array}{ll}
1, & \textrm{if $D_1 \ge \sqrt{D_2}$}, \\
\frac{\sqrt{D_2}}{\epsilon} \sim |z| \gtrsim \sqrt{|z|^2 + \frac{|z_l|^2}{\epsilon^2}} \gg 1, & \textrm{if $D_1 \le \sqrt{D_2}$}, \\
\end{array} \right.
\end{eqnarray*}
which proves our result. For the special case $|\alpha| = |\beta| = 1$, by \eqref{l1} and \eqref{index}, we have
\begin{align*}
|P_{\alpha,\iota}(z) | \lesssim d(z, t), \qquad \forall (z,t) \in \H(\mathcal{K}, \mathcal{A}),
\end{align*}
which yields $|\nabla p(g)| \lesssim d(g)$. Combining it with the result obtained above for $d(g)$ is large, we have shown $| \nabla \ln{p(g)} | \lesssim d(g)$.
~ \hspace*{20pt} ~ \hfill $\Box$

\medskip

\renewcommand{\theequation}{\thesection.\arabic{equation}}
\section{Sketch of proof for Grushin operators in two and three dimensions}  \label{st7}
\setcounter{equation}{0}

There are several places in \cite{Li12} where the assumption $n \ge 3$ is used: \\
(\romannumeral1) the estimates of $Q_2$ and $Q_{*2}$ in Section 3 and 4 respectively; \\
(\romannumeral2) the part ``Estimation de $|S_{m, 2}|$'' on p. 816; \\
(\romannumeral3) the Plancherel's theorem used in the $L^2(\R)$ sense in the second equality of (4.11), and the estimation (4.12) in Section 4.

Notice that (\romannumeral1) fails when $n = 1$ and $n =2$, and (\romannumeral2) and (\romannumeral3) fail only if $n = 1$. To overcome the difficulties mentioned above: a slight modification of \cite[Cas 2, p. 816]{Li12} suffices for (\romannumeral2);
for (\romannumeral1), it suffices to use the trick of integration by parts as in the treatment of $M_2$ and $M_2^*$ (see Subsections \ref{S52} and \ref{S52*} respectively); for (\romannumeral3), we have to establish \eqref{Plancherel} and \eqref{KE1} for $\nu = \frac{1}{2}$. In fact, we have the following more general result.

\begin{lemma}
\eqref{Plancherel} is valid for any $\Re \nu > 0$ and \eqref{KE1} is valid for any $\nu > 0$.
\end{lemma}

\noindent \textbf{Proof.}
For the first assertion, we observe that for fixed $r \geq 0$ and $b > 0$, the LHS and RHS of \eqref{Plancherel} are holomorphic in $\nu$ on $\{\nu : \Re \nu > 0 \}$ and they coincide for any $\nu \ge 1$, then the result follows from the uniqueness theorem of
holomorphic functions.

To prove the second assertion, we denote the RHS of \eqref{Plancherel} by $V = V(r,b;\nu)$ and we only need to prove for $\nu \in (0,1)$ and $\gamma_0 > 0$,
\begin{align*}
V \sim_{\gamma_0, \nu} \left(r + \sqrt{b}\right)^{\nu - 1}, \qquad \forall 0 \leq r \leq \gamma_0, \ 0 < b \leq \gamma_0.
\end{align*}

We restrict ourselves to the upper bound here since the lower bound follows from  a straightforward modification of \cite[``Minoration de $S_{m,1}$'' on p. 819]{Li12}.
First we recall (cf. \cite[\S 8.445, p. 919]{GR07}) that for fixed $\nu > 0$,
\begin{align*}
I_{\nu -1}(z) = \frac{1}{\Gamma(\nu)} \left(\frac{z}{2}\right)^{\nu - 1} (1 + o(1)), \qquad z \to 0, \ z \in \R,
\end{align*}
which implies that we have
\begin{align}\label{Gruest}
e^{-s} \left( \frac{s}{r} \right)^{\frac{\nu - 1}{2}} I_{\nu - 1}(2 \sqrt{r s})
\sim_{\nu} e^{-s} \left( \frac{s}{r} \right)^{\frac{\nu - 1}{2}} (\sqrt{r s})^{\nu - 1} = s^{\nu - 1}e^{-s}, \qquad \forall sr \leq 1.
\end{align}

Moreover, recall that (cf. \cite[\S8.451.5, p. 920]{GR07})
\begin{align*}
I_{\nu -1}(z) = e^{z} \frac{1}{\sqrt{2 \pi z}}(1 + o(1)),\qquad z \to \infty, \ z \in \R,
\end{align*}
which yields that for fixed $0 < \nu < 1$, we have for $s r > 1$ and $0 \leq r \leq \gamma_0$,
\begin{align*}
e^{-s} \left( \frac{s}{r} \right)^{\frac{\nu - 1}{2}} I_{\nu - 1}(2 \sqrt{r s})
&\sim_{\nu} e^{-s} \left( \frac{s}{r} \right)^{\frac{\nu - 1}{2}} e^{2 \sqrt{r s}} (\sqrt{r s})^{-\frac{1}{2}} \le s^{\nu - 1} e^{-s} (r s)^{\frac{1 - \nu}{2}} e^{2 \sqrt{r s}} \\
&\lesssim_{\nu, \gamma_0} s^{\nu - 1} \left(  s^{\frac{1 - \nu}{2}} e^{-s} e^{2 \sqrt{\gamma_0 s}} \right) \lesssim_{\nu, \gamma_0} s^{\nu - 1}.
\end{align*}
Combining it with \eqref{Gruest}, we obtain for fixed $0 < \nu < 1$ and $\gamma_0 > 0$,
\begin{align}\label{Gruest5}
e^{-s} \left( \frac{s}{r} \right)^{\frac{\nu - 1}{2}} I_{\nu - 1}(2 \sqrt{r s}) \lesssim_{\gamma_0,\nu} s^{\nu - 1}, \qquad \forall s > 0, \  0 \leq r \leq \gamma_0.
\end{align}

To derive the upper estimate, we use \eqref{Gruest5},
which yields
\begin{align*}
V \lesssim_{\gamma_0,\nu} \frac{1}{\sqrt{b}} \int_0^{+\infty} s^{\nu - 1} e^{- \frac{(s - r)^2}{4b}} \, ds
= \frac{1}{\sqrt{b}} \left( \int_0^{r + \sqrt{b}} + \int_{r + \sqrt{b}}^{+\infty}\right) = \Lambda_1 + \Lambda_2.
\end{align*}
The second term $\Lambda_2$ is easy to estimate:
\begin{align*}
\Lambda_2 \le \left(r + \sqrt{b}\right)^{\nu - 1} \frac{1}{\sqrt{b}} \int_{r + \sqrt{b}}^{+\infty}
 e^{- \frac{(s - r)^2}{4b}} \, ds \lesssim \left(r + \sqrt{b}\right)^{\nu - 1}.
\end{align*}
Now we are in the place to estimate $\Lambda_1$ and the proof is divided into two cases.\\
We first consider the case $r < 2\sqrt{b}$, in which we have
\begin{align*}
\Lambda_1 \le \frac{1}{\sqrt{b}} \int_0^{3\sqrt{b}} s^{\nu - 1} ds \sim_{\nu} (b)^{\frac{\nu - 1}{2} } \sim_{\nu} \left(r + \sqrt{b}\right)^{\nu - 1}.
\end{align*}
For the other case $r \ge 2\sqrt{b}$ we have
\begin{align*}
\Lambda_1 &\le \frac{1}{\sqrt{b}} \left(\int_0^{\frac{r}{2}} + \int_{\frac{r}{2}}^{r + \sqrt{b}}\right) \\
&\le \frac{1}{\sqrt{b}} e^{-\frac{r^2}{16b}} \int_0^{\frac{r}{2}} s^{\nu - 1} ds
+ \left(\frac{2}{r}\right)^{1 - \nu} \frac{1}{\sqrt{b}} \int_{\frac{r}{2}}^{r + \sqrt{b}}
e^{- \frac{(s - r)^2}{4b}} \, ds \\
&\lesssim_{\nu} (r)^{\nu - 1} \sim_{\nu} \left(r + \sqrt{b}\right)^{\nu - 1}.
\end{align*}
The results above complete the proof of the lemma.
~ \hspace*{20pt} ~ \hfill $\Box$

\begin{remark}
For the special case $\nu = \frac{1}{2}$, the lemma will be derived directly from the explicit expression of $I_{-\frac{1}{2}}$.
\end{remark}

\section*{Acknowledgement}
\setcounter{equation}{0}
This work is partially supported by
NSF of China (Grants  No. 11625102 and No. 11571077) and ``The Program of Shanghai Academic Research Leader'' (18XD1400700).
The authors would like to thank the anonymous referees for their many useful suggestions and valuable remarks which improve the writing of the paper, especially the constructive suggestion given by one referee which illuminates the results in Subsection \ref{futher}.

\mbox{}\\
Hong-Quan Li, Ye Zhang\\
School of Mathematical Sciences/Shanghai Center for Mathematical Sciences  \\
Fudan University \\
220 Handan Road  \\
Shanghai 200433  \\
People's Republic of China \\
E-Mail: hongquan\_li@fudan.edu.cn \quad or \quad hong\_quanli@yahoo.fr \\
17110180012@fudan.edu.cn \quad or \quad zhangye0217@126.com \mbox{}\\

\end{document}